\documentclass{amsart}

\title{The variance of closed geodesics in balls and annuli on the modular surface}

\author{Alexandre de Faveri}
\address{Department of Mathematics, Caltech, 1200 E. California Blvd., Pasadena, CA 91125, USA}
\email{\href{mailto:afaveri@caltech.edu}{afaveri@caltech.edu}}

\usepackage{amssymb, amsmath, amsthm}
\usepackage{hyperref, enumitem}
\usepackage[space]{cite}
\usepackage[centering]{geometry}
\usepackage{xcolor, graphicx}
\graphicspath{ {./images/} }
\usepackage{bbm}
\usepackage{mathrsfs}
\usepackage{pgf}
\usepackage{caption}
\usepackage{subcaption}
\usepackage[utf8]{inputenc}\DeclareUnicodeCharacter{2212}{-}

\newcommand{\dd}{\, d}
\newcommand{\R}{\mathbb{R}}
\renewcommand{\H}{\mathbb{H}}
\newcommand{\C}{\mathbb{C}}
\newcommand{\N}{\mathbb{N}}
\newcommand{\Z}{\mathbb{Z}}
\newcommand{\Q}{\mathbb{Q}}
\newcommand{\G}{\mathcal{G}}
\newcommand{\F}{\mathcal{F}}
\newcommand{\E}{\mathbb{E}}
\renewcommand{\P}{\mathbb{P}}
\newcommand{\eps}{\varepsilon}
\newcommand{\1}{\mathbbm{1}}
\newcommand{\norm}[1]{\left\lVert#1\right\rVert}

\newcommand{\dspec}{\, d_{\spec}}

\DeclareFontFamily{U}{mathx}{\hyphenchar\font45}
\DeclareFontShape{U}{mathx}{m}{n}{
      <5> <6> <7> <8> <9> <10>
      <10.95> <12> <14.4> <17.28> <20.74> <24.88>
      mathx10
      }{}
\DeclareSymbolFont{mathx}{U}{mathx}{m}{n}
\DeclareFontSubstitution{U}{mathx}{m}{n}
\DeclareMathAccent{\widecheck}{0}{mathx}{"71}
\DeclareMathAccent{\wideparen}{0}{mathx}{"75}

\newtheorem{theorem}{Theorem}
\newtheorem{lemma}{Lemma}

\newtheorem{corollary}{Corollary}

\newtheorem{remark}{Remark}

\newtheorem{conjecture}{Conjecture}

\DeclareMathOperator{\Cl}{Cl}
\DeclareMathOperator{\Var}{Var}
\DeclareMathOperator{\sym}{sym}
\DeclareMathOperator{\sgn}{sgn}
\DeclareMathOperator{\spec}{spec}

\DeclareMathOperator{\PSL}{PSL}
\DeclareMathOperator{\PSO}{PSO}
\DeclareMathOperator{\arccosh}{arccosh}
\DeclareMathOperator{\arcsinh}{arcsinh}
\DeclareMathOperator{\sech}{sech}

\begin{document}

\begin{abstract}
	We asymptotically estimate the variance for the distribution of closed geodesics in small random balls or annuli on the modular surface $\Gamma\backslash\H$. A probabilistic model in which closed geodesics are modeled using random geodesic segments is proposed, and we rigorously analyze this model using mixing of the geodesic flow in $\Gamma\backslash\H$. This leads to a conjecture for the asymptotic behavior of the variance, which unlike in previously explored cases is not equal to the expected value. We prove this conjecture for small balls and annuli, resolving a question left open by Humphries and Radziwi{\l}{\l}.
\end{abstract}

\maketitle


\section{Introduction}

Let $\Gamma := \PSL_2(\Z)$ denote the modular group and let $D > 0$ be a fundamental discriminant, meaning that $D$ is the discriminant of the real quadratic field $\Q(\sqrt{D})$. There is a well-known correspondence between narrow ideal classes in the narrow class group $\Cl_D^+$ of $\Q(\sqrt{D})$ and $\Gamma$-orbits of primitive irreducible integral binary quadratic forms $ax^2 + bxy + cy^2$ of discriminant $b^2 - 4ac = D$. Those, in turn, can also be associated to $\Gamma$-orbits of geodesics on the upper half-plane $\H$ with endpoints $\frac{-b \pm \sqrt{D}}{2a}$, or equivalently to the corresponding closed geodesics on the modular surface $\Gamma\backslash \H$. 

Denote the set of such closed geodesics of discriminant $D$ by $\Lambda_D$. Then $|\Lambda_D| = |\Cl_D^+| =: h_D^+$, and each closed geodesic in $\Lambda_D$ has length $2 \log{\eps_D^+}$, where $\eps_D^+ > 1$ is the smallest unit of positive norm in $\Q(\sqrt{D})$. The class number formula then gives
\begin{equation*}
    \sum_{\mathcal{C}\in \Lambda_D} \ell(\mathcal{C}) = h^+_D \cdot 2\log{\eps_D^+} = 2 \sqrt{D} L(1, \chi_D),
\end{equation*}
where $\chi_D$ is the primitive quadratic character modulo $D$ and $\ell(\mathcal{C}) := \int_\mathcal{C} ds$ denotes the length in $\H$, which is equipped with the hyperbolic metric and corresponding hyperbolic measure given respectively by
\begin{equation*}
    ds^2 := \frac{dx^2 + dy^2}{y^2} \quad \text{and} \quad d\mu(z) := \frac{dx \dd y}{y^2}
\end{equation*}
for $z = x+iy$. The bounds $D^{-\eps} \ll_\eps L(1, \chi_D) \ll \log{D}$ allow us to understand the total length quite well.

The elements of $\Lambda_D$ are expected to behave ``randomly'' in various senses (we will make this more precise below). In that direction, it is known that they become equidistributed in shrinking balls $B_R$: if we fix $\delta > 0$ and $w \in \Gamma \backslash\H$, then for $D^{-\frac{1}{18} + \delta} \ll R \ll 1$ we have
\begin{equation}\label{equidistribution_eq}
    \sum_{\mathcal{C}\in \Lambda_D} \ell(\mathcal{C} \cap B_{R}(w)) \sim \frac{\mu(B_R)}{\mu(\Gamma \backslash\H)} \sum_{\mathcal{C}\in \Lambda_D} \ell(\mathcal{C})
\end{equation}
as $D \to \infty$ through squarefree fundamental discriminants. Under the generalized Lindel\"of hypothesis we may replace the exponent $1/18$ by $1/6$, and equidistribution is expected to hold for exponents up to $1/2$. Such a result was first proved for fixed $R$ and with a congruence condition on $D$ by Skubenko \cite{Sku}, using Linnik's ergodic method \cite[Chapter VI]{Lin}. The congruence condition was only removed almost $30$ years later by Duke \cite{Duke_hyperb}, following a breakthrough of Iwaniec \cite{Iwa} (see \cite{ELMV} for a history of the problem). The result for shrinking $R$ mentioned above is given by Humphries \cite[Theorem 1.24]{Hum}, based on work of Young \cite{Young}. Analogous results are also available for geometric invariants in other contexts, such as Heegner points in $\H$ (corresponding to $D < 0$) and lattice points in spheres \cite{GF, DSP}, but we will restrict our attention to closed geodesics.

\begin{figure}[h]
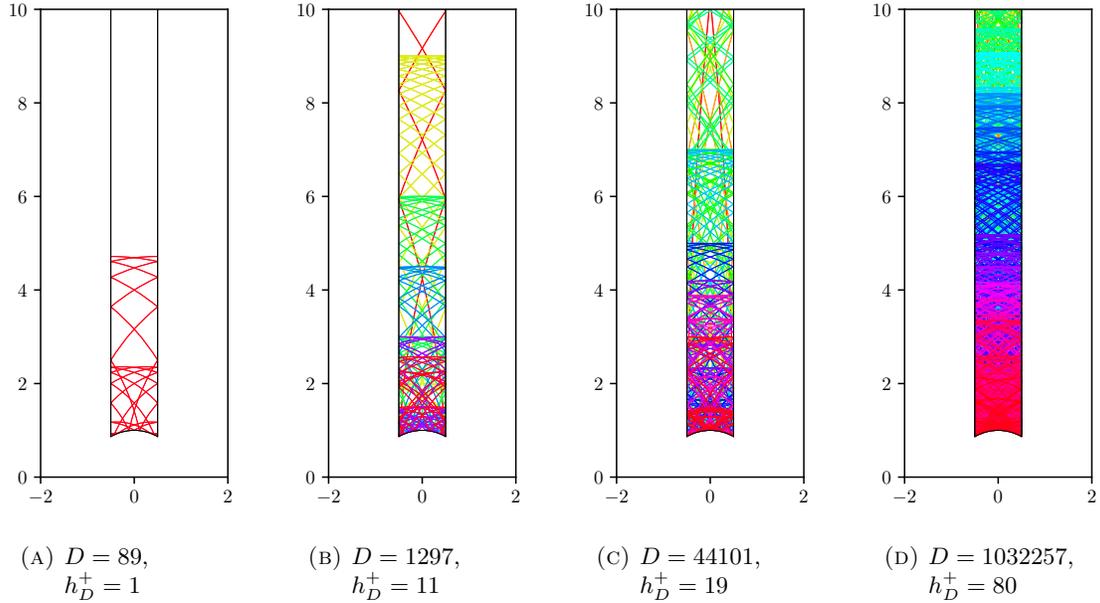

     \centering
     \begin{subfigure}[b]{0.24\textwidth}
         \centering
         \resizebox{\textwidth}{!}{\input{geodesics_89.pgf}}
         \captionsetup{format=hang, singlelinecheck=false, margin = 0.1\textwidth}
         \caption{$D = 89$,\newline$h^+_D = 1$}
     \end{subfigure}
     \hfill
     \begin{subfigure}[b]{0.24\textwidth}
         \centering
         \resizebox{\textwidth}{!}{\input{geodesics_1297.pgf}}
         \captionsetup{format=hang, singlelinecheck=false, margin = 0.1\textwidth}
         \caption{$D = 1297$,\newline$h^+_D = 11$}
     \end{subfigure}
     \hfill
     \begin{subfigure}[b]{0.24\textwidth}
         \centering
         \resizebox{\textwidth}{!}{\input{geodesics_44101.pgf}}
         \captionsetup{format=hang, singlelinecheck=false, margin = 0.1\textwidth}
         \caption{$D = 44101$,\newline$h^+_D = 19$}
     \end{subfigure}
     \hfill
     \begin{subfigure}[b]{0.24\textwidth}
         \centering
         \resizebox{\textwidth}{!}{\input{geodesics_1032257.pgf}}
         \captionsetup{format=hang, singlelinecheck=false, margin = 0.1\textwidth}
         \caption{$D = 1032257$,\newline$h^+_D = 80$}
     \end{subfigure}
        \caption{Closed geodesics in $\Lambda_D$}
\end{figure}

If one does not require equidistribution for every ball but instead is satisfied with a result covering almost all balls, then it is possible to go further. Considering a random variable given by the LHS of \eqref{equidistribution_eq}, where $w$ is distributed according to (a normalized version of) the measure $\mu$, it is tautological that the expected value is equal to the RHS of the same equation. One is then naturally led to consider the variance, which in the more general context of annuli $A_{r, R}(w)$ centered at $w \in \Gamma \backslash\H$, with inner radius $r$ and outer radius $R$, is given by
\begin{equation}\label{variance_eq}
    \Var(r, R; \Lambda_D) := \frac{1}{\mu(\Gamma\backslash\H)}\int_{\Gamma\backslash\H} \left( \sum_{\mathcal{C}\in \Lambda_D} \ell(\mathcal{C} \cap A_{r, R}(w)) - \frac{\mu(A_{r, R})}{\mu(\Gamma\backslash\H)} \sum_{\mathcal{C}\in \Lambda_D} \ell(\mathcal{C}) \right)^2 \dd\mu(w).
\end{equation}

Such an expression was first studied by Bourgain, Rudnick, and Sarnak \cite{BRS} in the context of lattice points in spheres. Based on probabilistic considerations, they conjectured that if the radii satisfy certain mild conditions, then the variance should be asymptotically equal to the corresponding expected value of the underlying random variable. An upper bound was then obtained assuming the generalized Lindel\"of hypothesis. 

Humphries and Radziwi{\l}{\l} \cite{HR19} were able to unconditionally prove the conjecture for certain very thin annuli, both in the case of lattice points in spheres and of Heegner points in $\H$. Furthermore, in the case of closed geodesics they obtained equidistribution for almost all annuli by showing that if $0 \leq r < R \ll 1$ and $D^{-1+\delta}\ll \mu(A_{r, R}) \ll 1$ for some fixed $\delta>0$, then for any fixed $c>0$,
\begin{equation*}
    \mu\left( \left\{ w \in \Gamma\backslash\H : \left| \frac{\mu(\Gamma\backslash\H)}{\mu(A_{r, R})} \frac{\sum_{\mathcal{C}\in \Lambda_D} \ell(\mathcal{C}\cap A_{r, R}(w))}{\sum_{\mathcal{C}\in \Lambda_D} \ell(\mathcal{C})} - 1 \right| > c \right\} \right) = o(1)
\end{equation*}
as $D \to \infty$ through squarefree fundamental discriminants. They did so by obtaining the bound $\Var(r, R; \Lambda_D) = o((\mu(A_{r, R}) \sqrt{D} L(1, \chi_D))^2)$ in this range, and indeed a careful examination of their method gives in particular
\begin{equation*}
    \Var(0, R; \Lambda_D) \ll_\eps  D^{\frac{1}{2}} R^{3-\eps}
\end{equation*}
for $R \ll D^{-\frac{5}{12}}$. This shows that the variance is not asymptotically equal to the expected value (which for balls is $D^{\frac{1}{2}+o(1)}R^2$, since $\mu(B_R) \asymp R^2$ for $R \ll 1$), as was the case for Heegner points in $\H$ and lattice points in spheres. A deviation of this kind is somewhat unexpected, since it implies better than ``square-root cancellation'' in \eqref{variance_eq}. However, in retrospect such a result is quite reasonable, since the geometric invariants have codimension $1$ in the case of closed geodesics, but $2$ in the other cases mentioned.

Given the discussion above, it is not completely clear what one should expect for the behavior of $\Var(r, R; \Lambda_D)$, and the purpose of this paper is to tackle this question. We start by proposing a probabilistic model, using geodesic segments of the appropriate length $2 \log{\eps^+_D}$ taken at random according to the Liouville measure in the unit tangent bundle of $\Gamma\backslash\H$, to model the elements of $\Lambda_D$ (see \autoref{random_section} for details). A rigorous analysis of this model turns out to be considerably more complicated than that for the geometric invariants of codimension $2$. We make critical use of a quantitative bound on the rate of mixing for the geodesic flow on the modular surface, combined with basic hyperbolic lattice point counting and some elementary hyperbolic geometry, to arrive at an asymptotic formula for the variance in the context of our probabilistic model. 

The main result in that direction is \autoref{random_thm}, where we show -- in the case of balls -- that for a single random geodesic segment of length $L$ in $\Gamma\backslash\H$, under mild conditions, the corresponding expression for the variance is $\sim \frac{16 L R^3}{\pi}$. For annuli, a certain special function $\mathbf{G}$ appears in the asymptotics (see \autoref{G_props_lemma} for its definition and key properties). We also refer to \autoref{heuristics_subsection} for a heuristic explanation of why the factor $R^3$ (instead of $R^2$) and the constant $\frac{16}{\pi}$ emerge in the asymptotics for this problem. Finally, it is worth pointing out that Luo and Sarnak \cite{LS04} have computed the quantum variance for the geodesic flow. In its classical incarnation, this variance is related to the spectral decomposition of our random model.

Using the analysis of the probabilistic model above, we are able to predict the asymptotic behavior of the variance for closed geodesics. In particular, in the case of balls we conjecture that if $0 < R \leq D^{-\delta}$ for some fixed $\delta > 0$, then
\begin{equation*}
    \Var(0, R; \Lambda_D) \sim \frac{64 \sqrt{D} L(1, \chi_D) R^3}{\pi}
\end{equation*}
as $D \to \infty$ through squarefree fundamental discriminants (see \autoref{main_conj} for the general case of annuli). Finally, our main result shows that the conjecture is true for balls of small radius.

\begin{corollary}\label{main_cor}
    Let $\delta > 0$ be given. If $0 < R \leq D^{-\frac{5}{12} - \delta}$, then as $D \to \infty$ through squarefree fundamental discriminants,
    \begin{equation*}
        \Var(0, R; \Lambda_D) \sim \frac{64 \sqrt{D} L(1, \chi_D) R^3}{\pi}.
    \end{equation*}
\end{corollary}

Indeed, \autoref{main_cor} is a particular case of \autoref{main_thm}, where we treat a wide class of annuli and the special function $\mathbf{G}$ appears, as expected. An interesting feature of the result is that the variance depends on the shape of the annulus, and not only on its area. The significance of the exponent $5/12$ and the obstacles towards extending the range of $R$ for which \autoref{main_cor} holds are discussed in \autoref{subconvexity_section}.

The proof of \autoref{main_thm} follows a completely different path than that of \autoref{random_thm}, and we instead apply the methods of \cite{HR19} to the case of closed geodesics. What allows us to prove a result for balls in this case is the presence of a different weight function than the one for Heegner points, due to the fact that the Gamma factors that arise when one expresses the relevant Weyl sums in terms of $L$-functions depend on the sign of $D$. The fact that the weight function decays faster is also a source of complications, since in our case the main contribution to the variance comes from forms with spectral parameter of size roughly between $1/R$ and $1/(R-r)$, as opposed to just around $1/(R-r)$ for Heegner points. This forces us to deal with the transition range $|x| \asymp 1$ for the Bessel function $J_0(x)$, where clear asymptotics are not available (see \autoref{bessel_rmk}). Thus instead of approximating with trigonometric functions, we carry the Bessel factors throughout the argument, and after certain integral transforms they are ultimately what gives rise to the special function $\mathbf{G}$ mentioned before in the asymptotics for the variance.

\subsection*{Acknowledgments}
I would like to thank my PhD advisor, Maksym Radziwi{\l}{\l}, for introducing me to this problem, and for general advice and encouragement. Thanks also to Valentin Blomer, Peter Humphries, Steve Lester, Carlos Matheus, and Ze{\'e}v Rudnick for helpful comments and suggestions on an earlier draft of this paper. I am also grateful to the two anonymous referees for their careful reading of the paper and useful comments, which greatly improved the readability of this work.


\section{Background and notation}

\subsection{Geometry of the upper half-plane}

The distance function $\rho:\H \times \H \to \R_{\geq0}$ and its more convenient proxy $u:\H \times \H \to \R_{\geq 0}$ are given by
\begin{equation*}
    \rho(z, w) := \log \left(\frac{|z- \overline{w}| + |z-w|}{|z- \overline{w}| - |z-w|}\right) \quad \text{and} \quad u(z, w) := \frac{|z-w|^2}{4 \Im(z) \Im(w)} = \sinh^2 \left(\frac{\rho(z, w)}{2}\right).
\end{equation*}
The group of isometries is $G:= \PSL_2(\R)$, which acts transitively through fractional linear transformations. The stabilizer of $i$ is $K:= \PSO_2(\R)$, so $g K \mapsto g i$ gives an identification $G/K \simeq \H$. 

Moreover, the corresponding action of $G$ on the unit tangent bundle $T^1(\H)$ (through the derivative map) is simply transitive, so if $v \in T^1(\H)$ denotes the unit tangent vector pointing up at $i$ then $g \mapsto g v$ gives an identification $G \simeq T^1(\H)$. More concretely, we can use the Iwasawa decomposition $G = N A K$, where
\begin{equation*}
    N := \left\{ \begin{pmatrix} 1 & t \\ 0 & 1 \end{pmatrix} : t \in \R \right\} \quad \text{and} \quad A := \left\{ \begin{pmatrix} a & 0 \\ 0 & a^{-1} \end{pmatrix}: a \in \R_{>0} \right\},
\end{equation*}
to describe this identification as 
\begin{equation*}
    \begin{pmatrix} 1 & x \\ 0 & 1 \end{pmatrix} \begin{pmatrix} y^{1/2} & 0 \\ 0 & y^{-1/2} \end{pmatrix} \begin{pmatrix} \cos(\frac{\theta}{2}) & \sin(\frac{\theta}{2}) \\ -\sin(\frac{\theta}{2}) & \cos(\frac{\theta}{2}) \end{pmatrix} \longleftrightarrow (z, \theta),
\end{equation*}
where $\theta$ is the angle with the unit tangent vector pointing up at $z = x+iy$. The derivative action of $G$ on $T^1(\H)$ becomes left multiplication in $G$ under the map described above, and the Liouville measure
\begin{equation*}
    d\nu(z, \theta) := \frac{dx \dd y}{y^2} \frac{d\theta}{2\pi}
\end{equation*}
on $T^1(\H)$ is invariant under this action of $G$, i.e. corresponds (up to a constant multiple) to the left-invariant Haar measure in $G$ under our identification. Furthermore since the group $G$ is unimodular, $\nu$ is also right-invariant.

\subsection{Geometry of the modular surface}

Let $X := \Gamma \backslash \H$ denote the modular surface, so that our previous identification quotients out to $X \simeq \Gamma \backslash G / K$, and similarly for the unit tangent bundle\footnote{Technically the modular surface has singularities at $i$ and $\frac{1+i \sqrt{3}}{2}$, since these points have nontrivial stabilizer in $\Gamma$. To correctly interpret the unit tangent bundle $T^1(X)$ we need to consider the orbifold structure of $X$, but this minor issue can be safely ignored for our purposes.} identification $T^1(X) \simeq \Gamma \backslash G$. The metric space structure of $X$ is obtained from the distance function
\begin{equation*}
    \widetilde{\rho}(\Gamma z, \Gamma w) := \min_{\gamma \in \Gamma} \rho(z, \gamma w).
\end{equation*}
Considering the usual (closure of a) fundamental domain
\begin{equation*}
    \F := \left\{z \in \H : |\Re(z)| \leq \frac{1}{2} \quad \text{and} \quad |z| \geq 1 \right\},
\end{equation*}
we can define measures $\widetilde{\mu}$ and $\widetilde{\nu}$ in $X$ and $T^1(X)$, respectively, by 
\begin{equation*}
    \widetilde{\mu}(\Gamma A) := \mu\left(\bigcup_{\gamma \in \Gamma} \gamma A \cap \F \right) \quad \text{and} \quad \widetilde{\nu}(\Gamma B) := \nu\left(\bigcup_{\gamma \in \Gamma} \gamma B \cap \pi^{-1}(\F) \right)
\end{equation*}
for measurable $A \subset \H$ and $B \subset T^1(\H)$, where $\pi:T^1(\H) \to \H$ is the projection map. In particular, $\widetilde{\nu}(T^1(X)) = \nu(\pi^{-1}(\F)) = \pi/3 = \mu(\F) = \widetilde{\mu}(X)$. Both measures are $G$-invariant under multiplication on the right, since the particular choice of fundamental domain turns out to be immaterial.

\subsection{Geodesic flow}

Given $t \in \R$, the geodesic flow $\G_t: T^1(\H) \to T^1(\H)$ is
\begin{equation*}
    \G_t(g) := g \begin{pmatrix} e^{t/2} & 0 \\ 0 & e^{-t/2} \end{pmatrix}
\end{equation*}
for $g \in G \simeq T^1(\H)$, and in geometric terms it amounts to parallel transport along the geodesic with starting point and direction given by the element of $T^1(\H)$ corresponding to $g$, for (hyperbolic) signed length $t$. The right-invariance of the Liouville measure $\nu$ implies that it is preserved by $\G_t$.

The geodesic flow clearly commutes with left multiplication by $G$ (and in particular by $\Gamma$), so it descends to a well-defined map $\widetilde{\G_t}: T^1(X) \to T^1(X)$ given by
\begin{equation*}
    \widetilde{\G_t}(\Gamma g) := \Gamma g \begin{pmatrix} e^{t/2} & 0 \\ 0 & e^{-t/2} \end{pmatrix}
\end{equation*}
for $\Gamma g \in T^1(X) \simeq \Gamma \backslash G$. Once again, $\widetilde{\G_t}$ preserves $\widetilde{\nu}$ and amounts to parallel transport by (hyperbolic) signed length $t$ along the corresponding geodesic in $X$.


\section{Estimates for the Selberg--Harish-Chandra transform}

\subsection{Definitions}

We follow \cite{HR19} with some minor modifications.

Let $k_{r, R}(u(z, w))$ be the identity function of the annulus 
\begin{equation*}
    A_{r, R}(w) := \left\{ z \in \H : r \leq \rho(z, w) \leq R\right\} = \left\{ z \in \H : \sinh^2 \left(\frac{r}{2}\right) \leq u(z, w) \leq \sinh^2 \left(\frac{R}{2}\right) \right\}
\end{equation*}
of hyperbolic volume
\begin{equation*}
    \mu(A_{r, R}) := \mu(A_{r, R}(w)) = 4 \pi \left( \sinh^2 \left(\frac{R}{2}\right) - \sinh^2 \left(\frac{r}{2}\right) \right),
\end{equation*}
that is,
\begin{equation*}
    k_{r, R}(t) := 
    \begin{cases}
        1 &\text{if } \sinh^2 \left(\frac{r}{2}\right) \leq t \leq \sinh^2 \left(\frac{R}{2}\right), \\
        0 &\text{otherwise}.
    \end{cases}
\end{equation*}
Observe that we use a different normalization from \cite{HR19} both here and in what follows below. Since $k_{r, R}(u(z, w))$ is a point-pair invariant, we can define the automorphic kernel $K_{r, R}: X \times X \to \R_{\geq 0}$ given by
\begin{equation*}
    K_{r, R}(z, w) := \sum_{\gamma\in\Gamma} k_{r, R}(u(z, \gamma w)).
\end{equation*}

The spectral expansion of this kernel involves the Selberg--Harish-Chandra transform $h_{r, R}$ of $k_{r, R}$, which is given by
\begin{equation}\label{SH_def_eq}
    h_{r, R}(t) := 2\pi \int_{0}^\infty P_{-\frac{1}{2} + it}(\cosh{\rho}) k_{r, R}\left(\sinh^2\left(\frac{\rho}{2}\right)\right) \sinh{\rho} \dd \rho = 2\pi \int_{r}^R P_{-\frac{1}{2} + it}(\cosh{\rho}) \sinh{\rho} \dd \rho, 
\end{equation}
where $P_\lambda$ is the Legendre function of the first kind.

\subsection{Bounds and asymptotics for \texorpdfstring{$h_{r, R}$}{h{r, R}}}

To understand the behavior of $h_{r, R}$ we express $P_{-\frac{1}{2} + it}$ in terms of Bessel functions, which will be more convenient to evaluate under the various integral transforms that will arise later.

\begin{lemma}[Hilb's formula {\cite[Lemma 2.24]{HR19}}]\label{Hilb_lemma}
    Fix $\eps > 0$. For $t \in \R$ and $0 < \rho < 1/\eps$,
    \begin{equation*}
        P_{-\frac{1}{2} + it}(\cosh{\rho}) = \sqrt{\frac{\rho}{\sinh{\rho}}} J_0(\rho t) + 
        \begin{cases}
            O(\rho^2) & \text{for } |t| \leq \frac{1}{\rho}, \\
            O_\eps\left(\frac{\sqrt{\rho}}{|t|^{3/2}}\right) &\text{for } |t| \geq \frac{1}{\rho} \geq \eps.
        \end{cases} 
    \end{equation*}
\end{lemma}

With this in mind, an asymptotic formula for $h_{r, R}$ easily follows. We restrict our attention to the case $R-r \gg R$, which will be relevant to us, but a similar statement also holds in the complementary case. 

\begin{lemma}\label{SH_asymp_lemma}
    Suppose that $0 \leq r < R \ll 1$ satisfy $R - r \gg R$, and $t \in \R$. Then 
    \begin{equation*}
        h_{r, R}(t) = 2 \pi \frac{R \cdot J_1(Rt) - r \cdot J_1(rt)}{t} + 
        \begin{cases}
            O(R^4) & \text{for } |t| \leq \frac{1}{R},\\
            O\left(\frac{R^{7/2}}{\sqrt{|t|}}\right) & \text{for } |t| \geq \frac{1}{R}.
        \end{cases}
    \end{equation*}
    Furthermore,
    \begin{equation*}
        h_{r, R}(t) \ll 
        \begin{cases}
            R^2 & \text{for } |t| \leq \frac{1}{R},\\
            \frac{\sqrt{R}}{|t|^{3/2}} & \text{for } |t| \geq \frac{1}{R}.
        \end{cases}
    \end{equation*}
\end{lemma}

\begin{proof}
Plugging \autoref{Hilb_lemma} into \eqref{SH_def_eq} gives
\begin{equation*}
    h_{r, R}(t) = 2\pi \int_{r}^R \sqrt{\rho \sinh{\rho}} \cdot J_0(\rho t) \dd \rho + 
    \begin{cases}
        O\left(R^4 \right) & \text{for } |t| \leq \frac{1}{R},\\
        O\left( \frac{R^{5/2}}{|t|^{3/2}} \right) & \text{for } |t| \geq \frac{1}{R}.
    \end{cases}
\end{equation*}
Using $\sinh{\rho} \ll \rho$, the bounds
\begin{equation}\label{J_0_bound_eq}
    J_0(x) = 
    \begin{cases}
        1 + O(x^2) & \text{for } |x| \leq 1,\\
        \sqrt{\frac{2}{\pi |x|}} \cos\left(|x| - \frac{\pi}{4}\right) + O\left(\frac{1}{|x|^{3/2}}\right) &\text{for } |x| \geq 1
    \end{cases}
\end{equation}
for $x \in \R$ \cite[8.411.1 and 8.451.1]{Integrals}, and integrating by parts in the case $|t| \geq 1/R$ (antidiferentiating the cosine term) gives the desired upper bound for $h_{r, R}$, as in \cite[Lemma 2.33]{HR19}. For the first asymptotic statement we use instead $\sinh{\rho} = \rho + O(\rho^3)$ combined with \eqref{J_0_bound_eq} to get
\begin{equation*}
    h_{r, R}(t) = 2 \pi \int_r^R \rho \cdot J_0(\rho t) \dd\rho + 
        \begin{cases}
            O\left(R^4 \right) & \text{for } |t| \leq \frac{1}{R},\\
            O\left( \frac{R^{7/2}}{\sqrt{|t|}} \right) & \text{for } |t| \geq \frac{1}{R}.
        \end{cases}
\end{equation*}
We can directly evaluate the remaining integral, since \cite[8.472.1]{Integrals} yields $(x J_1(x))' = x J_0(x)$, and the result follows.

\end{proof}

\begin{remark}\label{bessel_rmk}
    The reason we keep an expression with Bessel functions in the result above, instead of using \eqref{J_0_bound_eq} as in \cite[Lemma 2.27]{HR19} to write it in terms of simpler trigonometric functions, is that the main term in our variance computation will come roughly from $|t| \asymp 1/R$. This can be seen from the ranges of integration for the main term in \eqref{main_term_integral_range_eq}, as defined in \eqref{reference_eq}. In the case $R-r \ll R$ that range would be roughly $1/R \ll |t| \ll 1/(R-r)$, so either way we must deal with the transition range $|x| \asymp 1$ for $J_0(x)$, and \eqref{J_0_bound_eq} is not good enough to obtain asymptotics there. 
    
    In contrast, the main term in \cite[(7.18)]{HR19} -- with relevant ranges defined in \cite[(7.11)]{HR19} -- turns out to come roughly from $|t| \asymp 1/(R-r)$, which is much larger than $1/R$ (with the assumptions present there), so one still obtains an asymptotic for the Bessel function in the most important range. The main difference between the two cases is the presence of the extra weight $H(t)$ given by \eqref{gamma_factor_def_eq} in the spectral expansion of the variance for closed geodesics, which is not present in the case of Heegner points considered by Humphries and Radziwi{\l}{\l} (see \cite[Lemma 2.13]{HR19} for a comparison of the two weight functions).
\end{remark}


\section{Variance for random geodesic segments}\label{random_section}

Since the closed geodesics in $\Lambda_D$ are expected to behave in many aspects like ``random geodesics'', we will model them using uniformly distributed geodesic segments in $X$, so first we must understand the variance in that case.

By a geodesic segment of length $L$ in $X$, we mean a curve in $X$ of the form $\widetilde{\pi}\circ \widetilde{\G_t}(g)$ for $0 \leq t \leq L$ (observe that it is parametrized by hyperbolic arc length), where $g \in T^1(X)$ and $\widetilde{\pi}:T^1(X) \to X$ denotes the projection map. Uniform distribution means that the initial condition $g \in T^1(X)$ is distributed (up to normalization) according to the Liouville measure $\widetilde{\nu}$.

Given the discussion above, the random variable given by the length of the intersection between a random geodesic segment of length $L$ in $X$ with a random annulus $A_{r, R}$ in $X$ (with center distributed independently of the geodesic segment and according to the normalized measure in $X$) has variance
\begin{equation}\label{naive_variance_eq}
    \begin{split}
        \Var(r, R; L) :=& \int_{T^1(X)} \int_X \left(\int_0^L K_{r, R}(\widetilde{\pi}\circ \widetilde{\G_t}(g), w) \dd t - L \frac{\mu(A_{r, R})}{\widetilde{\mu}(X)}\right)^2 \frac{d\widetilde{\mu}(w)}{\widetilde{\mu}(X)} \frac{d\widetilde{\nu}(g)}{\widetilde{\nu}(T^1(X))} \\
        =& \int_{\pi^{-1}(\F)} \int_\F \left(\int_0^L \sum_{\gamma\in\Gamma} k_{r, R}( u(\pi\circ \G_t(g), \gamma w)) \dd t - L \frac{\mu(A_{r, R})}{\mu(\F)}\right)^2 \frac{d\mu(w)}{\mu(\F)} \frac{d\nu(g)}{\mu(\F)}.
    \end{split}
\end{equation}

\subsection{Heuristics}\label{heuristics_subsection} 

This discussion is partially inspired by the heuristics in \cite{Lalley}. 

For simplicity, we consider only the case of balls $B_R$ and assume $R = o(1)$. Let $\mathcal{S}$ be a random geodesic segment of length $L\gg 1$ on the modular surface $X$, and let $Y = Y(w) := \ell(\mathcal{S}\cap B_R(w))$, where $w$ is uniformly distributed in $X$. Also denote by $\mathcal{S}^R$ the ``tube'' of radius $R$ around $\mathcal{S}$. We wish to compute $\Var(Y)$, and will think of $\mathcal{S}$ as fixed but ``generic''. 

First suppose that $LR = o(1)$. Observe that $Y(w) = 0$ precisely for $w \not\in \mathcal{S}^R$, while $Y(w) \asymp R$ for most $w \in \mathcal{S}^R$ -- certainly for a typical $w \in \mathcal{S}^\frac{R}{2}$, since we expect this tube to have few self-intersections, as it has area $\asymp LR = o(1)$. Therefore, $\E(Y^2) \asymp LR^3$ and $\E(Y)^2 \asymp L^2R^4 = o(LR^3)$, so $\Var(Y) \asymp LR^3$. Furthermore, if for instance $\mathcal{S}^R$ has no self-intersections, then we can unfold it to $\H$ and obtain asymptotics for the variance using elementary hyperbolic geometry.

For the complementary case, suppose (say) that $LR \gg R^{\frac{1}{10}}$. We let $T \asymp R^{-\frac{9}{10}}$ and split $\mathcal{S}$ into $\frac{L}{T} \gg 1$ pieces $\mathcal{S}_1, \mathcal{S}_2, \dots, \mathcal{S}_{\frac{L}{T}}$ of length $T$. Since the geodesic flow is mixing (of all orders) and $T \to \infty$, we expect these segments to essentially behave independently. Let $Y_i = Y_i(w) := \ell(\mathcal{S}_i \cap B_R(w))$. Observe that $TR = o(1)$, so by the previous case we should have $\Var(Y_i) \asymp TR^3$. Then independence gives $\Var(Y) = \Var\left(\sum_{i=1}^\frac{L}{T} Y_i\right) \approx \sum_{i=1}^\frac{L}{T} \Var(Y_i) \asymp LR^3$, and it is reasonable to expect asymptotics for $\Var(Y)$ if we could obtain those for each $\Var(Y_i)$. 

In fact, one may heuristically determine the constant in $\Var(Y) \sim \frac{16 LR^3}{\pi}$ as follows. It suffices to consider the case $LR = o(1)$, by the argument using independence from the previous paragraph. The tube $\mathcal{S}^R$ has few self-intersections in that case, and the geometry of the problem is essentially Euclidean (as we work at the scale $R = o(1)$). Therefore, our situation can be modeled by the toy problem where the geodesic segment $\mathcal{S}$ is replaced by a straight line segment of length $L$ in $\R^2$, and the the point $w$ is randomized over some region of the appropriate area $\mu(X) = \frac{\pi}{3}$ that contains the (now Euclidean) tube $\mathcal{S}^R$. Since $\E(Y^2 \mid w \not\in \mathcal{S^R}) = 0$, we get $\E(Y^2) = \E(Y^2 \mid w \in \mathcal{S}^R) \cdot \P(w\in \mathcal{S}^R)$. Away from the endpoints of $\mathcal{S}$, $Y = 2\sqrt{R^2-x^2}$ depends only on the (signed) distance $x = x(w)$ from $w$ to the line that contains the segment $\mathcal{S}$, so we obtain
\begin{equation*}
    \E(Y^2 \mid w \in \mathcal{S}^R) \cdot \P(w\in \mathcal{S}^R) \approx \frac{1}{2R} \int_{-R}^{R} \left(2 \sqrt{R^2 - x^2}\right)^2 \dd x \cdot \frac{2RL}{\mu(X)} = \frac{8R^2}{3} \cdot \frac{6RL}{\pi} = \frac{16LR^3}{\pi}.
\end{equation*}

These heuristics provide a good intuition for the upcoming arguments in this section, but we will have to do something more complicated to effectively deal with self-intersections of $\mathcal{S}^R$.

\subsection{The cuspidal contribution}

Before delving into the variance computation, we need to make a small technical modification to \eqref{naive_variance_eq}, since with the current definition it turns out that $\Var(r, R;L) = \infty$ for $L \gg 1$. This is essentially due to the fact that the automorphic kernel $K_{r, R}(z, w)$ becomes quite large as $z$ and $w$ go towards the cusp together, so it is in particular not in $L^2(\F \times \F)$. This is the same issue that gives rise to continuous spectrum in the spectral resolution of the Laplacian in $X$.

For simplicity, consider the case of balls $B_R$, so $r = 0$. For $k \in \Z$ we have $\rho(w, w+k) \ll \sinh\left(\frac{\rho(w, w+k)}{2}\right) = \sqrt{u(w, w+k)} = \frac{|k|}{2\Im(w)}$, so for $\gg R \Im(w)$ values of $k \in \Z$ we have $\rho(w, w+k) \leq \frac{R}{4}$. Therefore, for all $g$ with $\pi(g) \in B_{\frac{R}{4}}(w)$ there are $\gg R \Im(w)$ values of $k\in\Z$ such that $\rho(\pi(g), w+k) \leq \rho(\pi(g), w) + \rho(w, w+k) \leq \frac{R}{2}$. Also observe that if $\rho(\pi(g), w + k) \leq \frac{R}{2}$ then $\pi(\G_t(g)) \in B_{R}(w+k)$ for all $0 \leq t \leq \frac{R}{2}$. 

We conclude that if $L \gg 1 \gg R$ then
\begin{equation*}
    \begin{split}
        \int_{\pi^{-1}(\F)} \int_{\F} \left( \int_{0}^L \sum_{\gamma\in\Gamma} k_{0, R}(u(\pi\circ\G_t(g), \gamma w)) \dd t \right)^2 \dd\mu(w) \dd\nu(g) \gg& \\
        \int_{\F} \int_{\F \cap B_{\frac{R}{4}}(w)} \left( \frac{R}{2} R \Im(w) \right)^2 \dd\mu(z) \dd\mu(w) \gg \int_{\frac{1}{R}}^\infty \frac{R}{y} (R^2 y)^2 \frac{d y}{y^2} =& \ \infty
    \end{split}
\end{equation*}
and this gives $\Var(0, R; T) = \infty$.

\subsection{The truncated variance}

In view of the necessity to exclude the contribution from the cusp, we let
$$\F_A := \{ z \in \F : \Im(z) \leq A\}$$
and consider averaging over annuli $A_{r, R}(w)$ only for $w\in \F_A$ instead of $w \in \F$, so that the relevant expression for the variance is
\begin{equation*}
    \Var_A(r, R; L) := \int_{\pi^{-1}(\F)} \int_{\F_A} \left(\int_0^L \sum_{\gamma\in\Gamma} k_{r, R}( u(\pi\circ \G_t(g), \gamma w)) \dd t - L \frac{\mu(A_{r, R})}{\mu(\F)}\right)^2 \frac{\dd\mu(w)}{\mu(\F_A)} \frac{\dd\nu(g)}{\mu(\F)}.
\end{equation*}

The asymptotic behavior of the expression above will involve a special function, so now we define it and express its key properties.

\begin{lemma}[Basic asymptotic properties of $\mathbf{G}$]\label{G_props_lemma}
    For $0 \leq w < 1$, let
    \begin{equation*}
        \mathbf{G}(w) := 1 + w^3 + (1-w^2)\mathbf{K}(w) -(1+w^2) \mathbf{E}(w),
    \end{equation*}
    where $\mathbf{K}$ and $\mathbf{E}$ are the complete elliptic integrals of the first and second kinds, respectively. Then 
    \begin{equation}\label{G_0_eq}
        \mathbf{G}(0) = 1,
    \end{equation}
    \begin{equation}\label{G_asymp_eq}
        \mathbf{G}(w) = \frac{3}{4} (1-w)^2 \log\left(\frac{2}{1-w}\right) + O\left((1-w)^2\right),
    \end{equation}
    \begin{equation}\label{G_bound_eq}
        \mathbf{G}(w) \gg (1-w)^2 \log\left(\frac{2}{1-w}\right),
    \end{equation}
    and
    \begin{equation}\label{G_der_bound_eq}
        \mathbf{G}'(w)  \ll (1-w) \log\left(\frac{2}{1-w}\right).
    \end{equation}
\end{lemma}

\begin{proof}
The definitions of $\mathbf{K}$ and $\mathbf{E}$ \cite[8.112]{Integrals} give \eqref{G_0_eq}, while \eqref{G_asymp_eq} and \eqref{G_der_bound_eq} follow from \cite[8.113.3 and 8.114.3]{Integrals}. Indeed, for $u := \sqrt{1-w^2}$, those give
\begin{equation*}
    \mathbf{K}(w) = \log\left(\frac{4}{u}\right) + \frac{u^2}{4}\log\left(\frac{4}{u}\right) + O(u^2)
\end{equation*}
and 
\begin{equation*}
    \mathbf{E}(w) = 1+\frac{u^2}{2}\log\left(\frac{4}{u}\right) -\frac{u^2}{4} + \frac{3u^4}{16}\log\left(\frac{4}{u}\right) + O(u^4),
\end{equation*}
so
\begin{equation*}
    \mathbf{G}(w) = 1+w^3 + u^2 \mathbf{K}(w) - (2-u^2)\mathbf{E}(w) = -1 + w^3 + \frac{3u^2}{2} + \frac{3u^4}{8}\log\left(\frac{4}{u}\right) + O(u^4).
\end{equation*}

Changing variables to $v := 1-w$, so $-1+w^3 = v^3+3v^2-3v$ and $u^2 = -v^2+2v$, we get
\begin{equation*}
    \begin{split}
        \mathbf{G}(w) &= v^3+3v^2-3v + \frac{3}{2}\left(-v^2+2v\right) + \frac{3}{8}\left(-v^2+2v\right)^2 \log\left(\frac{4}{\sqrt{v(1+w)}}\right) + O(v^2)\\
        & = \frac{3}{4}v^2 \log\left(\frac{1}{v}\right) + O(v^2),
    \end{split}
\end{equation*}
which gives \eqref{G_asymp_eq}. Similarly, the identity $\mathbf{G}'(w) = 3w(w - \mathbf{E}(w))$, which follows from \cite[8.123]{Integrals}, gives
\begin{equation*}
    \mathbf{G}'(w) \ll w- \mathbf{E}(w) \ll v \log\left(\frac{1}{v}\right) + O(v)
\end{equation*}
and we obtain \eqref{G_der_bound_eq}. This identity also shows that $\mathbf{G}'(w) \leq 0$, since $\mathbf{E}(w) \geq 1 > w$, so \eqref{G_bound_eq} follows from \eqref{G_asymp_eq} after choosing an appropriate cutoff.

\end{proof}

\begin{figure}[h]
    \centering
    \resizebox{0.75\textwidth}{!}{\input{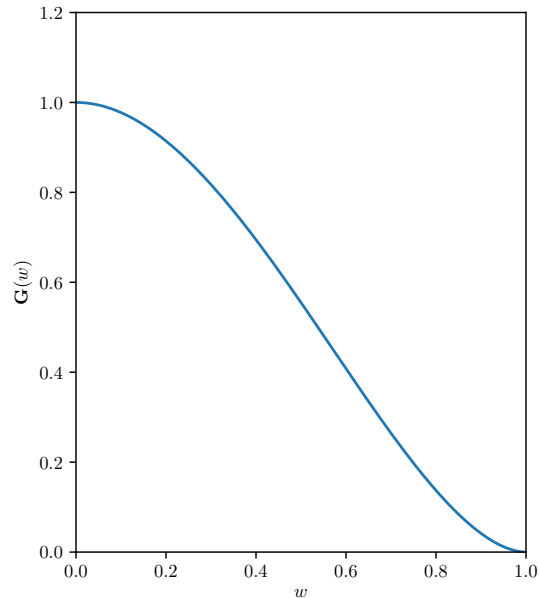}}
    \caption{Graph of $\mathbf{G}(w)$}
\end{figure}

\begin{remark}
    Perhaps the simplest way to understand the function $\mathbf{G}(w)$ geometrically is to describe it as follows: let $t \in [-1, 1]$ be uniformly distributed, and consider the (Euclidean) annulus $A'_{w, 1}(it) := \left\{ z \in \C : w \leq |z - it| \leq 1 \right\}$, for $w \in [0, 1]$. Let $Y_w$ be the (Euclidean) length of the intersection of $A'_{w, 1}(it)$ with the real axis. Then an explicit computation shows that
    \begin{equation*}
        \mathbf{G}(w) = \frac{\E(Y_w^2)}{\E(Y_0^2)} = \frac{3}{8}\E(Y_w^2) = \frac{3}{8}\cdot\frac{1}{2}\int_{-1}^{1} \left(2\sqrt{1-t^2} - \1_{|t|\leq w} \cdot 2\sqrt{w^2-t^2}\right)^2 \dd t.
    \end{equation*}
\end{remark}

We are ready to state the main result of this section, omitting the dependence on a parameter $t$ that governs all the asymptotic statements below (meaning that the quantities $L, r, R, A$ are all functions of $t$, and asymptotic notations such as $o(\dots)$ or $\sim$ should be interpreted in the limit as $t \to \infty$).

\begin{theorem}\label{random_thm} 
    Suppose that $L \gg 1$, $0 \leq r < R = o(1)$, and $1 \ll \log{A} = o\left(R^{-1} \log^{-1} \left({\frac{1}{R-r}}\right)\right)$, so in particular we require $\log \left(\frac{1}{R-r}\right) = o\left(\frac{1}{R}\right)$. Then
    \begin{equation*}
        \Var_A(r, R; L) \sim \frac{16 L R^3}{\pi} \mathbf{G}\left(\frac{r}{R}\right).
    \end{equation*}
    In particular, for balls we get
    \begin{equation*}
        \Var_A(0, R; L) \sim \frac{16 L R^3}{\pi},
    \end{equation*}
    and for thin annuli (i.e. such that $R-r = o(R)$) satisfying the restrictions above we get
    \begin{equation*}
        \Var_A(r, R; L) \sim \frac{12 L R (R-r)^2}{\pi} \log\left(\frac{R}{R-r}\right).
    \end{equation*}
\end{theorem}

\subsection{Auxiliary results}

An important ingredient for \autoref{random_thm} will be the fact that the geodesic flow is mixing, and in fact it is so with an exponential rate, due to a theorem of M. Ratner \cite{Ratner}. We will use the following effective version of Ratner's result, due to C. Matheus and adapted here to the modular group $\Gamma$.

\begin{lemma}[Exponential mixing for the geodesic flow {\cite[Corollary 2.1]{Matheus}}]\label{mixing_lemma}
    Let $\phi, \psi \in L^2(X)$ be such that $\int_X \phi \dd\mu = \int_X \psi \dd\mu = 0$. Then
    \begin{equation*}
        \int_{T^1(X)} \phi(\widetilde{\pi}(g))\cdot \psi(\widetilde{\pi} \circ \widetilde{\G_t}(g)) \dd\widetilde{\nu}(g) \ll \norm{\phi}_{L^2(X)}\cdot \norm{\psi}_{L^2(X)}\cdot (|t| + 1) e^{-\frac{|t|}{2}}.
    \end{equation*}
\end{lemma}

In order to deal with the problem of self-intersections alluded to in our heuristic discussion, we will need two basic observations regarding the distribution of orbits of $\Gamma$ acting on $\H$. The first one is the following standard density estimate.

\begin{lemma}[Density of hyperbolic lattice points {\cite[Lemma 2.11]{Iwaniec}}]\label{lattice_density_lemma} 
    If $z \in \H$, $w \in \F$, and $\delta > 0$ then
    \begin{equation*}
        \left|\{ \gamma \in \Gamma : u(z, \gamma w) \leq \delta \}\right| \ll \sqrt{\delta(\delta + 1)} \Im(w) + 1.
    \end{equation*}
\end{lemma}

The second observation about the orbits of $\Gamma$ in $\H$ deals with the minimum spacing between distinct points in such an orbit. It formalizes the idea that the spacing can only be small if either the orbit comes close to a point of $\H$ with nontrivial stabilizer in $\Gamma$, or if it has a point very high up towards the cusp.

\begin{lemma}[Minimum spacing of hyperbolic lattice points]\label{min_lemma}
    If $w \in \F$ then
    \begin{equation*}
        \min_{1 \not= \gamma\in \Gamma} \rho(w, \gamma w) \gg \min\left\{\rho(w, i), \rho(w, j), \rho(w, j'), \frac{1}{\Im(w)}\right\}
    \end{equation*}
    where $j := \frac{1 + i\sqrt{3}}{2}$ and $j' := \frac{-1 + i\sqrt{3}}{2}$.
\end{lemma}

\begin{proof}
Since the minimum is $\leq \rho(w, w+1) \ll 1$, it suffices to show that
\begin{equation*}
    \min_{1 \not= \gamma\in \Gamma} u(w, \gamma w) \gg \min\left\{u(w, i), u(w, j), u(w, j'), \frac{1}{\Im(w)^2}\right\}.
\end{equation*}
Write $w = x+iy$, $1 \not= \gamma = \begin{pmatrix} a & b \\ c & d \end{pmatrix} \in \Gamma$ for an element that attains the minimum, and $U^2 := u(w, \gamma w)$ with $1 \gg U > 0$.

\begin{itemize}[leftmargin=*, align=left]
    \item[\textit{Case 1:}] $y \geq 2$.
    
        If $c \not= 0$ we have $\Im(\gamma w) = \frac{y}{|cw+d|^2} \leq \frac{y}{(cy)^2}\leq \frac{1}{2}$, which gives $U \gg 1$. If $c = 0$ we have $\gamma w = w \pm b$ and $b \not= 0$, which gives $U \gg \frac{1}{\Im(w)}$. In any case, the result holds for $y \geq 2$.
    \item[\textit{Case 2:}] $y < 2$.
    
        Observe that since $w \in \F$ we have $y \geq \frac{\sqrt{3}}{2}$. If $|c| \geq 2$ then as before $\Im(\gamma w) = \frac{y}{|cw +d|^2} \leq \frac{1}{c^2 y} \leq \frac{1}{2 \sqrt{3}}$, so $U \gg 1$. If $c = 0$ we also get $U \gg 1$ as before, so $|c| = 1$ and we can assume that $c = 1$ since the entries of $\gamma$ are only determined up to flipping all the signs.
        
        Now, $\gamma = \begin{pmatrix} a & ad-1 \\ 1 & d \end{pmatrix}$ and
        \begin{equation*}
            U^2 = \frac{|w - \gamma w|^2}{4 y \Im(\gamma w)} = \frac{|w - \gamma w|^2 |w+d|^2}{4y^2} = \frac{|w(w+d) - (aw + ad - 1)|^2}{4y^2} \gg |z|^2
        \end{equation*}
        for $z := w(w+d) - (aw + ad - 1)$, so that
        \begin{equation*}
            \Re(z) = x^2 - y^2 + (d-a)x - ad + 1 \quad \quad \text{and} \quad \quad \Im(z) = 2xy + (d-a)y.
        \end{equation*}
        
        If $|d-a| \geq 2$ then since $|x| \leq \frac{1}{2}$, as $w\in \F$, we get $|\Im(z)| \geq y \gg 1$ and therefore $U \gg 1$.
        
        If $d = a \not= 0$ then we have $\Re(z) = x^2 - y^2 - a^2 + 1 \leq \frac{1}{4} - \frac{3}{4} = -\frac{1}{2}$, so $U \gg 1$. 
        
        If $d = a = 0$ then $U^2 \gg \Im(z)^2 \gg x^2$ and $U^2 \gg \Re(z)^2 = (x^2 - y^2 + 1)^2$, which gives $|x| \ll U$ and $y^2 = 1 + x^2 + O(U) = 1 + O(U)$ since $U \ll 1$, so that $|y-1| \ll U$. We conclude that $u(w, i) \asymp x^2 + (y-1)^2 \ll U^2$, as desired.
        
        If $d - a = 1$ then we can assume that $d = 0$ or $1$, otherwise $ad \geq 2$ and $\Re(z) = x^2 - y^2 + x -ad + 1 \leq \frac{1}{4} - \frac{3}{4} + \frac{1}{2} - 2 + 1 = -1$, so $U \gg 1$. For $d=0$ we get $U \gg |\Im(z)| \gg \left|x+\frac{1}{2}\right|$, so $x = -\frac{1}{2} + O(U)$ and then looking at the real part we get $y^2 = 1 + x + x^2 + O(U) = \frac{3}{4} + O(U)$, once again since $U \ll 1$. This gives $y = \frac{\sqrt{3}}{2} + O(U)$ and therefore $u(w, j) \asymp \left(x+\frac{1}{2}\right)^2 + \left(y-\frac{\sqrt{3}}{2}\right)^2 \ll U^2$, as desired. For $d = 1$ the exact same reasoning shows that $u(w, j) \ll U^2$.
        
        Finally, if $d-a = -1$ then an argument analogous to the previous paragraph, but exchanging $x$ with $-x$, gives $u(w, j') \ll U^2$, so we have covered all possibilities and the result follows.
\end{itemize}
\end{proof}

The last ingredients necessary to prove \autoref{random_thm} are bounds and asymptotics for integral expressions that measure the lengths of intersections between geodesics and annuli in $\H$, averaged over various parameters. We deal with those geometric quantities in the next two lemmas, and emphasize that the results are analogous (except for large distances) to those for the Euclidean version of the problem, in which straight lines intersect Euclidean annuli.

\begin{lemma}[Average intersection of geodesics through $z$ with $A_{r, R}(w)$]\label{geodesic_int_bound_lemma}
    For any $z, w \in \H$, and $0 \leq r < R \ll 1$, if we denote $D := \rho(z, w)$ then
    \begin{equation*}
        \Theta_{r, R}(z, w) := \int_{0}^{2\pi} \int_{0}^\infty k_{r, R}(u(\pi\circ \G_t(z, \theta), w)) \dd t \dd \theta \ll
        \begin{cases}
            (R-r) \log\left(\frac{2R}{R-r}\right) & \text{if } D < 2R, \\
            \frac{R(R-r)}{D} & \text{if } 2R \leq D < 1, \\
            \frac{R(R-r)}{e^{D}} & \text{if } 1 \leq D.
        \end{cases}
    \end{equation*}
\end{lemma}

\begin{proof}
Since the geodesic flow is parametrized by arc length, the system of coordinates $(t, \theta)$ corresponds to geodesic polar coordinates centered at $z$, therefore the hyperbolic measure becomes $d\mu = 2 \sinh{t} \dd t \dd\theta$. 

If $D \geq 2R$, then observing that the integrand is simply the indicator function of the annulus $A_{r, R}(w)$ and that $\sinh{t} \geq \sinh(D-R) \gg \sinh{D}$ for all points $(t, \theta)$ inside it (since by the triangle inequality $t = \rho(z, (t, \theta)) \geq \rho(z, w) - \rho((t, \theta), w) \geq D-R$), we get
\begin{equation*}
    \Theta_{r, R}(z, w) \ll \int_{0}^{2\pi} \int_{0}^\infty k_{r, R}(u(\pi\circ \G_t(z, \theta), w)) \frac{\sinh{t}}{\sinh{D}} \dd t \dd \theta \ll \frac{\mu(A_{r, R})}{\sinh{D}}
\end{equation*}
and the result follows. If $R-r \gg R$ then the result for $D < 2 R$ is trivial since the integral over $t$ is always $\leq 2R$ by the triangle inequality. Therefore we can assume that $r \geq \frac{R}{2}$, and then a slight modification of the argument above also takes care of $D \leq \frac{R}{4}$, since $\sinh{t} \gg \sinh{R}$ in that case.

We are left with the trickiest case $r \geq \frac{R}{2}$ and $\frac{R}{4} \leq D \leq 2R$. The issue here is that the intersection of each geodesic with the annulus $A_{r, R}(w)$ no longer has length $\ll R-r$ when the latter is thin -- in fact the length can be $\gg \sqrt{R(R-r)}$. In what follows it is worth keeping in mind that since $R \ll 1$ the geometry is roughly Euclidean.

First let us change variables, shifting $\theta$ so that it corresponds to the angle with the geodesic from $z$ to $w$ (instead of with the vertical line). Since $D \ll R$, we can choose a sufficiently small (absolute) $\eps > 0$ such that if $|\sin{\theta}| < \eps$ then each $\theta$ contributes $\ll R-r$. Indeed, the integrand for each $\theta$ is now the length of the intersection of the one-sided geodesic determined by $(z, \theta)$ with the annulus $A_{r, R}(w)$. If $d$ is the (orthogonal) distance between that geodesic and $w$, then the hyperbolic law of sines gives $\sinh{d} = \sinh{D} \sin{\theta}$, which implies $d \ll R \sin{\theta}$. The length of the intersection of the corresponding two-sided geodesic with $A_{r, R}(w)$ is $2(\ell_R - \ell_r)$, where we define $2 \ell_R$ as the length of the intersection of that two-sided geodesic with $B_R(w)$, and similarly for $2\ell_r$ (both intersections will be non-empty for sufficiently small $\eps>0$, as we assume $D \asymp R \asymp r$). By the hyperbolic law of cosines we have (see \autoref{length_fig})
\begin{equation*}
    \cosh{\ell_R} = \frac{\cosh{R}}{\cosh{d}} \quad \quad \text{and} \quad \quad \cosh{\ell_r} = \frac{\cosh{r}}{\cosh{d}}.
\end{equation*}

\begin{figure}[h]
    \centering
    \resizebox{0.9\textwidth}{!}{\input{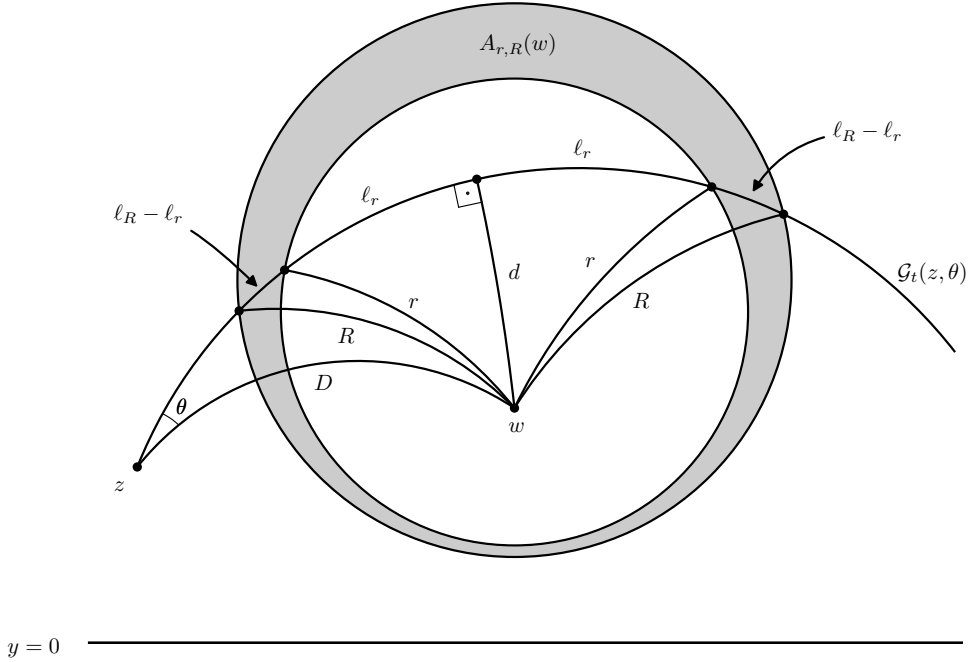}}
    \caption{Lengths in intersection of $\mathcal{G}_t(z, \theta)$ with $A_{r, R}(w)$}
    \label{length_fig}
\end{figure}

Choosing $\eps>0$ sufficiently small so that $d \leq \frac{R}{4} \leq \frac{r}{2}$ we get $\ell_r \gg R$ and therefore by the MVT
\begin{equation*}
    \ell_R - \ell_r \ll \frac{\cosh{\ell_R} - \cosh{\ell_r}}{\sinh{\ell_r}} \ll \frac{\cosh{R} - \cosh{r}}{R\cosh{d}} \ll R-r.
\end{equation*}

Now, for the remaining angles $\theta$ satisfying $|\sin{\theta}| \geq \eps$ we will fix the radius $t$ and evaluate the angular contribution. The values of $t$ for which $\partial B_{t}(z)$ intersects only one of $\partial B_R(w)$ or $\partial B_{r}(w)$ contribute $\ll R-r$ (bounding the integral over $\theta$ trivially), so they may be excluded and we can assume that both circles are intersected. Let $0 \leq \theta_r < \theta_R \leq \pi$ be the angles corresponding to the intersection points in the upper half, so the integral over $\theta$ with no restrictions contributes $2(\theta_R - \theta_r)$ and by the hyperbolic law of cosines we have (see \autoref{angles_fig})
\begin{equation*}
    \cos{\theta_R} = \frac{\cosh{t}\cosh{D} - \cosh{R}}{\sinh{t}\sinh{D}} \quad \quad \text{and} \quad\quad \cos{\theta_r} = \frac{\cosh{t}\cosh{D} - \cosh{r}}{\sinh{t}\sinh{D}}. 
\end{equation*}

\begin{figure}[h]
    \centering
    \resizebox{0.9\textwidth}{!}{\input{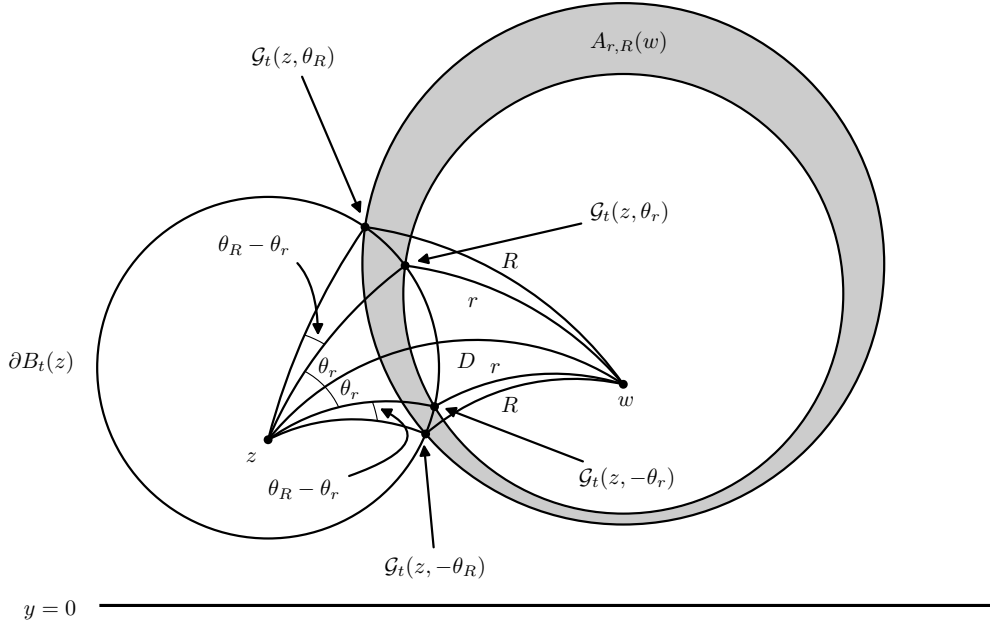}}
    \caption{Angles in intersection of $\partial B_t(z)$ with $A_{r, R}(w)$}
    \label{angles_fig}
\end{figure}

If $[\theta'_r, \theta'_R] := [\theta_r, \theta_R] \cap [\arcsin{\eps}, \pi - \arcsin{\eps}]$, then since we have already excluded angles $\theta$ with $|\sin{\theta}| < \eps$ by bounding the corresponding terms as in the previous part of the argument, the contribution of the remaining terms corresponding to $t$ is actually just $2(\theta'_R - \theta'_r)$. But since the sine of both angles is bounded away from zero we can use the MVT to get
\begin{equation*}
    \theta'_R - \theta'_r \ll \frac{\cos{\theta'_r} - \cos{\theta'_R}}{\sin{\eps}} \ll \cos{\theta_r} - \cos{\theta_R} = \frac{\cosh{R} - \cosh{r}}{\sinh{t} \sinh{D}} \ll \frac{R-r}{t},
\end{equation*}
as $D \gg R$. 

Finally, we integrate over $t \leq D+R \ll R$, where we can assume that $t \geq R-r$ since those radii trivially contribute $\ll R-r$. In conclusion,
\begin{equation*}
    \Theta_{r, R}(z, w) \ll R-r + \int_{R-r}^{D+R} \frac{R-r}{t} \dd t \ll (R-r) \log\left(\frac{2R}{R-r}\right),
\end{equation*}
as desired.

\end{proof}

\begin{lemma}[Main term computation]\label{main_geodesic_int_lemma}
    For any $g \in T^1(\H)$ and $0 \leq r < R = o(1)$,
    \begin{equation*}
        \int_{\H} k_{r, R}(u(\pi(g), w)) \left( \int_{-\infty}^{\infty} k_{r, R}(u(\pi \circ \G_t(g), w)) \dd t \right)  \dd\mu(w) \sim \frac{16 R^3}{3} \mathbf{G}\left(\frac{r}{R}\right).
    \end{equation*}
\end{lemma}

\begin{proof}
Observe that the LHS is independent of $g$, since the integral over $w$ is invariant under isometries, so denoting the whole expression by $I_{r, R}$ we can average over the geodesic segment of length $S > 0$ (which we will choose to be sufficiently large later) starting at $(i, 0) \in T^1(\H)$ to get
\begin{equation*}
    I_{r, R} = \frac{1}{S} \int_{\H} \left( \int_{0}^S k_{r, R}(u(i e^s, w)) \dd s \right) \left( \int_{-\infty}^{\infty} k_{r, R}(u(i e^t, w)) \dd t \right)  \dd\mu(w).
\end{equation*}
Given $D> 0$, let 
\begin{equation*}
    A(D, S) := \left\{r e^{i \delta} \in \H :  1\leq r \leq e^S \quad \text{and} \quad \frac{\pi}{2} - \theta_D \leq \delta \leq \frac{\pi}{2} + \theta_D \right\},
\end{equation*}
for $0 \leq \theta_D \leq \frac{\pi}{2}$ defined by $\sin{\theta_D} = \tanh{D}$ (see \autoref{region_fig}). It will be important to compute $\mu(A(D, S))$ for the computation of $I_{r, R}$ that follows below, so we do that now and come back to the integral afterwards. 

The locus of points $z\in \H$ with (orthogonal) distance to the line $y = 0$ equal to $D$ is the pair of straight half-lines through the origin with angles $\theta_D$ and $-\theta_D$ with the vertical. The (hyperbolic) arc length parametrization of the half-line corresponding to $\theta_D$ is $z(s) = e^{s \cos{\theta_D}} e^{i \left(\frac{\pi}{2} - \theta_D\right)}$. A computation shows that its geodesic curvature is constant equal to $\sin{\theta_D} = \tanh{D}$ (see for instance the discussion after \cite[Corollary 4]{Cza}, where our situation corresponds to $a = 1$ and $\alpha = \frac{\pi}{2}-\theta_D$). The region $A(D, S)$ has as boundaries two geodesics (Euclidean circles with center at the origin) and the two straight lines through the origin with angles $\theta_D$ and $-\theta_D$ with the vertical line $y=0$. Examining the arc length parametrization we see that each of those has length $\frac{S}{\cos{\theta_D}}$, so the total geodesic curvature along the boundary of $A(D, S)$ in the positive direction is 
\begin{equation*}
    2 S \tan{\theta_D} = 2 S \frac{\tanh{D}}{\sqrt{1 - \tanh^2{D}}} = 2 S \sinh{D}.
\end{equation*}

\begin{figure}[h]
    \centering
    \resizebox{0.9\textwidth}{!}{\input{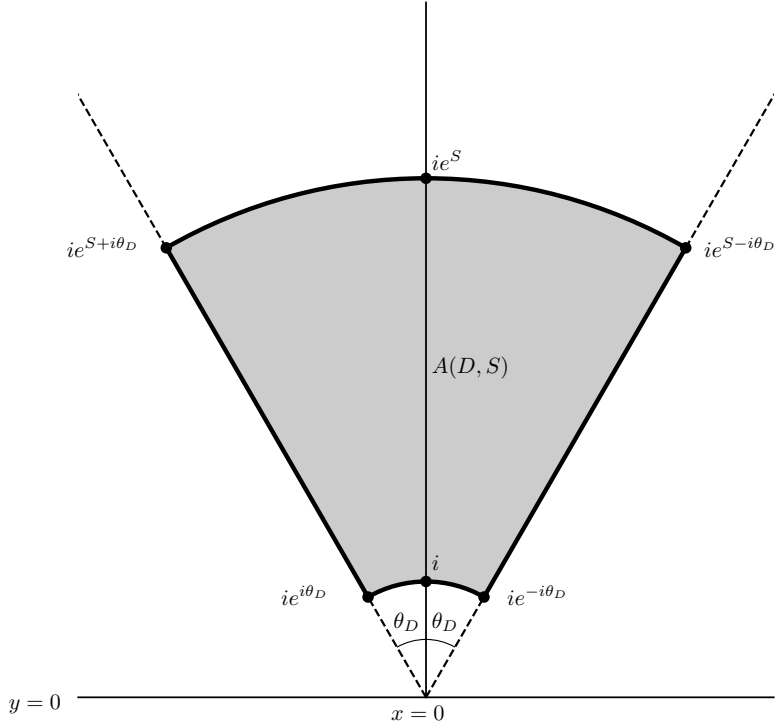}}
    \caption{The region $A(D, S)$}
    \label{region_fig}
\end{figure}

Since the four external angles of $A(D, S)$ are equal to $\frac{\pi}{2}$, denoting by $K=-1$ the Gaussian curvature of $A(D, S)$ and by $k_g$ the geodesic curvature of $\partial A(D,S)$, the Gauss-Bonnet theorem gives
\begin{equation*}
    \mu(A(D, S)) = -\int_{A(D, S)} K \dd\mu = -2\pi + 4 \frac{\pi}{2} + \int_{\partial A(D, S)} k_g \dd s = 2 S \sinh{D}.
\end{equation*}

Observe that $\int_{0}^S k_{r, R}(u(i e^s, w)) \dd s$ is nonzero only for $w \in A(R, S) \cup B_R(i) \cup B_R(i e^S)$, since this is the locus of points within distance $\leq R$ from the geodesic segment between the points $i$ and $ie^S$. Furthermore, if $w \not\in B_R(i) \cup B_R(i e^S)$ then $\int_{0}^S k_{r, R}(u(i e^s, w)) \dd s$ is equal to $\int_{-\infty}^\infty k_{r, R}(u(i e^s, w)) \dd s$, which is the length of the intersection of the vertical line $y = 0$ with the annulus $A_{r, R}(w)$ (and therefore trivially $\ll R$). As discussed in the proof of \autoref{geodesic_int_bound_lemma}, the hyperbolic law of cosines shows that for points at distance $D$ from the vertical this length is $2 (\ell_R - \ell_r)$ if $0 \leq D \leq r$ and $2 \ell_R$ for $r \leq D \leq R$, where 
\begin{equation*}
    \cosh{\ell_R} = \frac{\cosh{R}}{\cosh{D}} \quad \quad \text{and} \quad \quad \cosh{\ell_r} = \frac{\cosh{r}}{\cosh{D}}.
\end{equation*}

We conclude that
\begin{equation*}
    I_{r, R} = \frac{1}{S}\int_{A(R, S)} \left( \int_{-\infty}^\infty k_{r, R}(u(i e^t, w)) \dd t \right) ^2 \dd\mu(w) + O\left(\frac{R^4}{S}\right).
\end{equation*}
Choosing say $S \gg (R-r)^{-2}$ we see that the error term is $o(R(R-r)^2)$ and will be negligible. The remaining term can be written as
\begin{equation}\label{arccos_integral_eq}
    \begin{split}
        \frac{4}{S} &\int_{0}^R \left(\arccosh\left( \frac{\cosh{R}}{\cosh{D}} \right) - \1_{D \leq r} \cdot \arccosh\left( \frac{\cosh{r}}{\cosh{D}} \right)\right)^2  \dd \left(\mu(A(D, S) )\right) \\
        = 8 &\int_{0}^R \left(\arccosh\left( \frac{\cosh{R}}{\cosh{D}} \right) - \1_{D \leq r} \cdot \arccosh\left( \frac{\cosh{r}}{\cosh{D}} \right)\right)^2 \cosh{D} \dd D \\
        \sim 8 &\int_{0}^r \log^2\left( \frac{\cosh{R} + \sqrt{\cosh^2{R} - \cosh^2{D}}}{\cosh{r} + \sqrt{\cosh^2{r} - \cosh^2{D}}} \right) \dd D \\ 
        & \quad \quad \quad \quad \quad \quad \quad \quad \quad \quad \quad \quad + 8 \int_{r}^R \log^2\left( \frac{\cosh{R} + \sqrt{\cosh^2{R} - \cosh^2{D}}}{\cosh{D}}\right) \dd D
    \end{split}
\end{equation}
using $\arccosh{x} = \log\left(x + \sqrt{x^2-1} \right)$ and $\cosh{D} = 1 + O(R^2)$. 

Denoting 
\begin{equation*}
    f_{r, R}(D) := \sqrt{\cosh^2{R} - \cosh^2{D}} - \sqrt{\cosh^2{r} - \cosh^2{D}} = \sqrt{\sinh^2{R} - \sinh^2{D}} - \sqrt{\sinh^2{r} - \sinh^2{D}}
\end{equation*}
for $0 \leq D \leq r$, we see that it is increasing and therefore the MVT gives
\begin{equation}\label{basic_sqrt_bound_eq}
    R-r \ll f_{r, R}(D) \ll \sqrt{R (R-r)}.
\end{equation}
Similarly, $\cosh{r} + \sqrt{\cosh^2{r} - \cosh^2{D}} = 1 + O(R)$ and $\cosh{R} - \cosh{r} = O(R(R-r))$, so
\begin{equation*}
    \frac{\cosh{R} + \sqrt{\cosh^2{R} - \cosh^2{D}}}{\cosh{r} + \sqrt{\cosh^2{r} - \cosh^2{D}}} = 1 + f_{r, R}(D) (1 + O(R)) + O(R(R-r))
\end{equation*}
and therefore with the aid of \eqref{basic_sqrt_bound_eq} we obtain
\begin{equation*}
    \log^2\left( \frac{\cosh{R} + \sqrt{\cosh^2{R} - \cosh^2{D}}}{\cosh{r} + \sqrt{\cosh^2{r} - \cosh^2{D}}} \right) = f_{r, R}(D)^2 (1 + O(R)) + O(R^2(R-r)^2).
\end{equation*}
The same sort of analysis for $g_R(D) := \sqrt{\cosh^2{R} - \cosh^2{D}} \ll \sqrt{R (R-r)}$ when $r \leq D \leq R$, separating into cases depending on whether $g_R(D)$ is larger than $R-r$ or not, gives
\begin{equation*}
    \log^2\left( \frac{\cosh{R} + \sqrt{\cosh^2{R} - \cosh^2{D}}}{\cosh{D}} \right) = g_{R}(D)^2 (1 + O(R)) + O(R(R-r)^2).
\end{equation*}

Plugging those into \eqref{arccos_integral_eq} we get an error term $\ll R^2 (R-r)^2 = o(R(R-r)^2)$ and an integral term 
\begin{equation*}
    \begin{split}
        &\sim 8\int_0^r f_{r, R}(D)^2 \dd D + 8\int_r^R g_R(D)^2 \dd D \\
        &\sim 8 \int_{0}^{\sinh{r}} f_{r, R}(\arcsinh{x})^2 \dd x + 8 \int_{\sinh{r}}^{\sinh{R}} g_R(\arcsinh{x})^2 \dd x,
    \end{split}
\end{equation*}
changing variables to $x = \sinh{D}$ and observing that $\cosh{D} = 1 + O(R^2)$. The result can be written, by \cite[3.155.8]{Integrals}, as
\begin{equation*}
    \begin{split}
        &\frac{16}{3} \left(\sinh^3{R} + \sinh^3{r}\right) - 16 \int_0^{\sinh{r}} \sqrt{(\sinh^2{R} - x^2) (\sinh^2{r} - x^2)} \dd x \\
        = \  &\frac{16}{3} \left(\sinh^3{R} + \sinh^3{r}\right)- \frac{16}{3} \sinh{R} \left( \left(\sinh^2{R} + \sinh^2{r}\right) \mathbf{E}\left(\frac{\sinh{r}}{\sinh{R}}\right) - \left(\sinh^2{R} - \sinh^2{r}\right) \mathbf{K}\left( \frac{\sinh{r}}{\sinh{R}} \right) \right) \\
        = \ & \frac{16 \sinh^3{R}}{3} \mathbf{G}\left(\frac{\sinh{r}}{\sinh{R}}\right) \sim \frac{16 R^3}{3} \mathbf{G}\left(\frac{\sinh{r}}{\sinh{R}}\right),
    \end{split}
\end{equation*}
where we use the notation of \autoref{G_props_lemma}. A computation with Taylor series gives
\begin{equation*}
    \frac{\sinh{r}}{\sinh{R}} = \frac{r}{R} + O\left(R(R-r)\right),
\end{equation*}
so \eqref{G_der_bound_eq} and the MVT give
\begin{equation*}
    \mathbf{G}\left(\frac{\sinh{r}}{\sinh{R}}\right) = \mathbf{G}\left(\frac{r}{R}\right) + O\left((R-r)^2 \log\left( \frac{2R}{R-r} \right) \right).
\end{equation*}

We conclude that 
\begin{equation*}
    I_{r, R} \sim \frac{16 R^3}{3} \mathbf{G}\left(\frac{r}{R}\right) + o(R(R-r)^2) \sim \frac{16 R^3}{3} \mathbf{G}\left(\frac{r}{R}\right)
\end{equation*}
by \eqref{G_bound_eq}, as desired.

\end{proof}

\subsection{Putting it all together: the proof of \autoref{random_thm}}

\begin{proof}[Proof of \autoref{random_thm}]
By absolute convergence, we can freely exchange the order of integration and write
\begin{equation*} 
    \Var_A(r, R; L) = \int_{\F_A} \int_{T^1(X)} \left(\int_{0}^L K_{r, R}(\widetilde{\pi}\circ \widetilde{\G}_{t}(g), w) \dd t - L\frac{\mu(A_{r, R})}{\mu(\F)}\right)^2 \frac{d\widetilde{\nu}(g)}{\mu(\F)} \frac{d\mu(w)}{\mu(\F_A)}.
\end{equation*}
The inner integral over $g \in T^1(X)$ is, after changing variables, equal to
\begin{equation}\label{cross_term_eq}
    \int_0^L \int_0^L \int_{T^1(X)} \phi(\widetilde{\pi}(g)) \cdot \phi(\widetilde{\pi} \circ \widetilde{\G}_{t-t'}(g)) \frac{d\widetilde{\nu}(g)}{\mu(\F)} \dd t' \dd t
\end{equation}
for $\phi \in L^2(X)$ given by $\phi(z) := K_{r, R}(z, w) - \frac{\mu(A_{r, R})}{\mu(\F)}$, so that $\int_X \phi \dd\mu = 0$ and
\begin{equation}\label{norm_bound}
    \begin{split}
        \norm{\phi}^2_{L^2(X)} = & \int_{\F} K_{r, R}^2(z, w) \dd\mu(z) - \left(\frac{\mu(A_{r, R})}{\mu(\F)}\right)^2 < \int_{\F} K_{r, R}^2(z, w) \dd\mu(z)\\
        \leq & \max_{\xi \in \F} \left|\left\{\gamma\in\Gamma: u(\xi, \gamma w) \leq \sinh^2 \left(\frac{R}{2}\right) \right\}\right| \cdot \int_{\F} K_{r, R}(z, w) \dd\mu(z) \\
        \ll & (R \Im(w) + 1) \frac{\mu(A_{r, R})}{\mu(\F)} \ll (R \Im(w) + 1) R(R-r)
    \end{split}
\end{equation}
by \autoref{lattice_density_lemma}. Let $L \geq T \gg 1$. If $|t - t'| > T$ we can use \autoref{mixing_lemma} in the inner integral of \eqref{cross_term_eq}, inputting the bound \eqref{norm_bound}, to conclude that the contribution of all such terms to $\Var_A(r, R; L)$ is
\begin{equation}\label{error_1_eq}
    \begin{split}
        &\ll LR(R-r) \int_{\F_A} (R \Im(w) + 1) \dd\mu(w) \int_T^\infty (x + 1) e^{-\frac{x}{2}} \dd x  \\
        &\ll L R(R-r) T e^{-\frac{T}{2}} (R \log{A} + 1) \ll L R(R-r) e^{-\frac{T}{4}}
    \end{split}
\end{equation}
if $T < L$, and it is $=0$ if we choose $T = L$ (since no such terms exist in that case).

The remaining set of $|t - t'| \leq T$ has measure $\ll L T$, so replacing $\phi(z)$ with $K_{r, R}(z, w)$ in \eqref{cross_term_eq} we pick up an error term
\begin{equation}\label{error_2_eq}
    \ll L T R^2(R-r)^2,
\end{equation}
and what remains is
\begin{equation*}
    \begin{split}
        \int_{0}^L \int_{\max \{t'-T, 0\}}^{\min \{t'+T, L\}} \int_{T^1(X)} K_{r, R}(\widetilde{\pi}(g), w) \cdot K_{r, R}(\widetilde{\pi}\circ \widetilde{\G}_{t-t'}(g), w) \frac{d\widetilde{\nu}(g)}{\mu(\F)} \dd t \dd t' \\
        =  \int_{T^1(X)} K_{r, R}(\widetilde{\pi}(g), w) \left(\int_{-T}^{T} (L - |t|) \cdot K_{r, R}(\widetilde{\pi}\circ \widetilde{\G}_{t}(g), w)\dd t \right) \frac{d\widetilde{\nu}(g)}{\mu(\F)}.
    \end{split}
\end{equation*}
Inserting the integral above into the expression for the variance and expanding the automorphic kernel, we are left with
\begin{equation}\label{interm_var_eq}
    \begin{split}
         \int_{\F_A} & \int_{\pi^{-1}(\F)} \sum_{\gamma' \in \Gamma} k_{r, R} (u(\pi(g), \gamma' w))  \left(\int_{-T}^T (L-|t|) \sum_{\gamma \in \Gamma}  k_{r, R} (u(\pi\circ \G_t(g), \gamma w)) \dd t \right) \frac{d\nu(g)}{\mu(\F)} \frac{d\mu(w)}{\mu(\F_A)} \\
        &= \int_{\F_A} \sum_{\gamma' \in \Gamma} \int_{\pi^{-1}(\F \cap A_{r, R}(\gamma' w))} \sum_{\gamma \in \Gamma} \int_{-T}^T (L-|t|) \cdot k_{r, R} (u(\pi\circ \G_t(g), \gamma w)) \dd t \frac{d\nu(g)}{\mu(\F)} \frac{d\mu(w)}{\mu(\F_A)}.
    \end{split}
\end{equation}
Now, let $\mathcal{M}$ denote the terms corresponding to $\gamma = \gamma'$, which is where the main contribution will come from, and let $\mathcal{E}$ denote all other terms. 

Decomposing $g = (z, \theta)$ for $z \in \F\cap A_{r, R}(\gamma'w)$ and $0 \leq \theta < 2\pi$, what is left in \eqref{interm_var_eq} corresponding to the terms in $\mathcal{E}$ is
\begin{equation}\label{E_bound_eq}
    \ll L \int_{\F_{A}} \sum_{\gamma'\in \Gamma}\int_{\F \cap A_{r, R}(\gamma'w)} \int_0^{2\pi} \sum_{\gamma'\not= \gamma\in \Gamma}\int_{-T}^T k_{r, R}(u(\pi\circ \G_t(z, \theta), \gamma w)) \dd t \dd\theta \dd\mu(z)\dd\mu(w).
\end{equation}
For given $w \in \F_{A}$, $\gamma' \in \Gamma$, and $z \in \F\cap A_{r, R}(\gamma' w)$, we can use the notation of \autoref{geodesic_int_bound_lemma} to bound
\begin{equation}\label{angle_bound_eq}
    \int_0^{2\pi} \sum_{\gamma'\not= \gamma \in \Gamma} \int_{-T}^{T} k_{r, R}(u(\pi\circ \G_t(z, \theta), \gamma w)) \dd t \dd \theta \ll \sum_{\gamma'\not= \gamma\in \Gamma} \1_{\rho(z, \gamma w) < T+R} \cdot \Theta_{r, R}(z, \gamma w).
\end{equation}
Denote
\begin{equation*}
    h(\gamma', w, z; D) := \left|\left\{ \gamma' \not= \gamma \in \Gamma : \rho(z, \gamma w) \leq D \right\}\right|, \quad \quad f(w; D) := \left| \left\{ 1 \not= \gamma \in \Gamma : \rho(w, \gamma w) \leq D \right\}\right|.
\end{equation*}
The fact that $z \in A_{r, R}(\gamma' w)$ implies $h(\gamma', w, z; D) \leq f(w; D+R)$. The next step is to apply \autoref{geodesic_int_bound_lemma} in the equation below, where in order to simplify the notation we keep a term $(R-r)^{-1} \log^{-1}\left(\frac{2R}{R-r}\right)$ in the left and adjust the bounds in each range so that the boundary terms cancel out after integration by parts. This gives
\begin{equation}\label{angle_integral_eq}
    \begin{split}
        &(R-r)^{-1} \log^{-1}\left(\frac{2R}{R-r}\right) \sum_{\gamma'\not= \gamma\in \Gamma} \1_{\rho(z, \gamma w) < T+R} \cdot \Theta_{r, R}(z, \gamma w) \\ 
        \ll\ &\int_{0^-}^{R} \dd h(\gamma', w, z; D) + \int_{R}^{1} \frac{R}{D} \dd h(\gamma', w, z; D) + \int_{1}^{T+R} \frac{R}{e^{D-1}} \dd h(\gamma', w, z; D) \\
        =\ &\int_{R}^1 \frac{R}{D^2} h(\gamma', w, z; D) \dd D + \int_{1}^{T+R} \frac{R}{e^{D-1}} h(\gamma', w, z; D)  \dd D + \frac{R}{e^{T+R-1}} h(\gamma', w, z; T+R) \\
        \ll\ &\int_{R}^1 \frac{R}{D^2} f(w; D + R) \dd D + \int_{1}^{T+R} \frac{R}{e^{D}} f(w; D+R)  \dd D + \frac{R}{e^{T}} f(w; T+2R).
    \end{split}
\end{equation}

We now denote 
\begin{equation*}
    m(w) := \min_{1 \not= \gamma\in \Gamma} \rho(w, \gamma w) = \min \{ D \in \R : f(w; D) > 0 \}
\end{equation*}
and consider two different cases for $w \in \F_A$:
\begin{itemize}[leftmargin=*, align=left]
    \item[\textit{Case 1:}] $w \in \F_A \cap B_{1/100}(q)$ for some $q \in \{i, j, j'\}$.
        
        In this case \autoref{lattice_density_lemma} and \autoref{min_lemma} give respectively
        \begin{equation*}
            f(w; D) \ll 
            \begin{cases}
                1 & \text{for } D \leq 1, \\
                e^D & \text{for } D > 1,
            \end{cases}
            \quad \text{and} \quad m(w) \gg \rho(w, q),
        \end{equation*}
        so that \eqref{angle_integral_eq} is
        \begin{equation*}
            \ll \int_{\max\{m(w), R\}}^1 \frac{R}{D^2} \dd D + \int_1^{T+R} R \dd D + R \ll  \frac{R}{\max\{\rho(w, q), R\}} + TR.
        \end{equation*}
        Plugging this into \eqref{angle_bound_eq} and then \eqref{E_bound_eq}, we see that the total contribution to $\mathcal{E}$ in this case is 
        \begin{equation}\label{E_bound_1_eq}
            \begin{split}
                &\ll L (R-r) \log\left(\frac{2R}{R-r}\right) \int_{\F_A \cap B_{1/100}(q)} \sum_{\gamma' \in \Gamma} \int_{\F \cap A_{r, R}(\gamma'w)} \left(\frac{R}{\max\{\rho(w, q), R\}} + TR\right) \dd\mu(z) \dd\mu(w) \\
                &\ll L (R-r) \log\left(\frac{2R}{R-r}\right) \int_{B_{1/100}(q)} \left(\frac{R}{\max\{\rho(w, q), R\}} + TR\right) \int_\F K_{r, R}(z, w) \dd\mu(z) \dd\mu(w) \\
                &\ll LR(R-r)^2 \log\left(\frac{2R}{R-r}\right) (R^2 + R + TR) \ll L T R^2 (R-r)^2 \log\left(\frac{2R}{R-r}\right).
            \end{split}
        \end{equation}
    
    \item[\textit{Case 2:}] $w \in \F_A$ but $w \not\in B_{1/100}(q)$ for any $q \in \{i, j, j'\}$.
        
        In this case \autoref{lattice_density_lemma} and \autoref{min_lemma} give respectively
        \begin{equation*}
            f(w; D) \ll 
            \begin{cases}
                D \Im(w) + 1 & \text{for } D \leq 1, \\
                e^D \Im(w) & \text{for } D > 1,
            \end{cases}
            \quad \text{and} \quad m(w) \gg \frac{1}{\Im(w)},
        \end{equation*}
        so that \eqref{angle_integral_eq} is
        \begin{equation*}
            \begin{split}
                &\ll \int_{\max\{m(w), R\}}^1 \frac{R}{D^2} (D\Im(w) + 1) \dd D + \int_1^{T+R} R \Im(w) \dd D + R \Im(w) \\
                &\ll R \Im(w) \log\left(\frac{1}{R}\right) + RT \Im(w).
            \end{split}
        \end{equation*}
        Denote $\F^*_A := \F_A \setminus \bigcup_{q \in \left\{i, j, j'\right\}} B_{1/100}(q)$. Plugging the bound above into \eqref{angle_bound_eq} and then \eqref{E_bound_eq}, we see that the total contribution to $\mathcal{E}$ in this case is 
        \begin{equation}\label{E_bound_2_eq}
            \begin{split}
                &\ll L (R-r) \log\left(\frac{2R}{R-r}\right) \int_{\F^*_A} \sum_{\gamma' \in \Gamma} \int_{\F \cap A_{r, R}(\gamma'w)} \left(R \Im(w) \log\left(\frac{1}{R}\right) + RT \Im(w) \right) \dd\mu(z) \dd\mu(w) \\
                &\ll L (R-r) \log\left(\frac{2R}{R-r}\right) \int_{\F_A} \left(R \Im(w) \log\left(\frac{1}{R}\right) + RT \Im(w)\right) \int_\F K_{r, R}(z, w) \dd\mu(z) \dd\mu(w) \\
                &\ll L R^2 (R-r)^2 \log\left(\frac{2R}{R-r}\right) \left(\log\left(\frac{1}{R}\right) + T \right) \log{A}.
            \end{split}
        \end{equation}
\end{itemize}

Collecting the error terms \eqref{error_1_eq} and \eqref{error_2_eq}, and the estimates for $\mathcal{E}$ in \eqref{E_bound_1_eq} and \eqref{E_bound_2_eq}, we conclude that 
\begin{equation*}
    \begin{split}
        \Var_A(r, R; L) &= \mathcal{M} + \mathcal{E} + O\left(L T R^2(R-r)^2 + \1_{T \not= L} \cdot L R(R-r) e^{-\frac{T}{4}}\right) \\
        &= \mathcal{M} + O\left(L R^2 (R-r)^2 \log\left(\frac{2R}{R-r}\right) \left(\log\left(\frac{1}{R}\right) + T \right) \log{A} + \1_{T \not= L} \cdot L R(R-r) e^{-\frac{T}{4}}\right).
    \end{split}
\end{equation*}
Choosing $T = \min\left\{8 \log\left(\frac{1}{R-r}\right), L\right\}$, recalling that $\log{A} = o\left(R^{-1} \log^{-1} \left(\frac{1}{R-r}\right) \right)$ we see that the error term is $o\left(L R (R-r)^2 \log\left(\frac{2R}{R-r}\right)\right)$.

We are left with computing the main term 
\begin{equation*}
    \mathcal{M} := \int_{\F_A} \sum_{\gamma' \in \Gamma} \int_{\pi^{-1}(\F \cap A_{r, R}(\gamma' w))} \int_{-T}^{T} (L-|t|) \cdot k_{r, R}(u(\pi \circ \G_t(g), \gamma' w)) \dd t \frac{d\nu(g)}{\mu(\F)} \frac{d\mu(w)}{\mu(\F_A)}.
\end{equation*}
Observe that since $\pi(g) \in A_{r, R}(\gamma' w)$ we can restrict the integral over $t$ to $[-2R, 2R]$, as the intersection of a geodesic with an annulus $A_{r, R}$ in $\H$ is contained in a segment of the geodesic of length $\leq 2R$. Then $L - |t| = L - O(R) \sim L$. Therefore, 
\begin{equation*}
    \mathcal{M} \sim L \int_{\F_A} \sum_{\gamma' \in \Gamma} \int_{\pi^{-1}(\F)} k_{r, R}(u(\pi(g), \gamma'w)) \left(\int_{-2R}^{2R} k_{r, R}(u(\pi \circ \G_t(g), \gamma' w)) \dd t \right)  \frac{d\nu(g)}{\mu(\F)} \frac{d\mu(w)}{\mu(\F_A)}.
\end{equation*}

We have $u(\pi \circ \G_t(g), \gamma' w) = u(\gamma'^{-1} \pi \circ \G_t(g), w) = u(\pi \circ \G_t(\gamma'^{-1}g), w)$, and the measure $\nu$ is (left) $G$-invariant, so it is possible to unfold the integral over $g \in \pi^{-1}(\F)$ to get
\begin{equation*}
    \mathcal{M} \sim \frac{L}{\mu(\F)\mu(\F_A)} \int_{\F_A} \int_{G} k_{r, R}(u(\pi(g), w)) \left(\int_{-2R}^{2R} k_{r, R}(u(\pi \circ \G_t(g), w)) \dd t \right)  \dd\nu(g) \dd\mu(w).
\end{equation*}
By the same argument via $G$-invariance, the integral over $g \in G$ is independent of $w$, so we can replace the domain $\F_A$ with $\F$ in the display above, multiplying by $\frac{\mu(\F_A)}{\mu(\F)}$, to get
\begin{equation*}
    \begin{split}
        \frac{L}{\mu(\F)^2} \int_{\F} \int_{G} k_{r, R}(u(\pi(g), w)) \left(\int_{-2R}^{2R} k_{r, R}(u(\pi \circ \G_t(g), w)) \dd t \right) \dd\nu(g)\dd\mu(w).
    \end{split}
\end{equation*}
Now we revert the unfolding process in the integral over $g\in G$, obtaining
\begin{equation*}
    \mathcal{M} \sim \frac{L}{\mu(\F)^2} \int_{\F} \sum_{\gamma \in \Gamma} \int_{\pi^{-1}(\F)} k_{r, R}(u(\pi(g), \gamma w)) \left(\int_{-2R}^{2R} k_{r, R}(u(\pi \circ \G_t(g), \gamma  w)) \dd t \right)   \dd\nu(g)\dd\mu(w).
\end{equation*}

Next we can complete the inner integral to $t \in (-\infty, \infty)$ with no loss, as was previously discussed, and after unfolding (this time the integral over $w \in \F$) we are left with
\begin{equation*}
    \mathcal{M} \sim \frac{L}{\mu(\F)^2} \int_{\pi^{-1}(\F)} \int_{\H} k_{r, R}(u(\pi(g), w)) \left( \int_{-\infty}^{\infty} k_{r, R}(u(\pi \circ \G_t(g), w)) \dd t \right) \dd\mu(w) \dd\nu(g). 
\end{equation*}

Applying \autoref{main_geodesic_int_lemma} and \eqref{G_bound_eq}, we conclude that
\begin{equation*}
    \mathcal{M} \sim \frac{16 L R^3}{3 \mu(\F)} \mathbf{G}\left(\frac{r}{R}\right) = \frac{16 L R^3}{\pi} \mathbf{G}\left(\frac{r}{R}\right) \gg L R (R-r)^2 \log\left(\frac{2R}{R-r}\right).
\end{equation*}
Therefore,
\begin{equation*}
    \Var_A(r, R; L) = \mathcal{M} + o\left(L R (R-r)^2 \log\left(\frac{2R}{R-r}\right)\right) \sim \frac{16 L R^3}{\pi} \mathbf{G}\left(\frac{r}{R}\right),
\end{equation*}
as desired. Combining this with \eqref{G_0_eq} and \eqref{G_asymp_eq} finishes the proof of the theorem.
    
\end{proof}

\begin{remark}\label{mixing_rmk}
    It may be possible to remove the technical condition $\log\left(\frac{1}{R-r}\right) = o\left(\frac{1}{R}\right)$ in \autoref{random_thm}, extending the result to all $L \gg 1$, $0\leq r < R = o(1)$, and $1 \ll \log{A} = o\left(R^{-1} \log^{-1}\left(\frac{1}{R}\right)\right)$. That is because we use a somewhat simplified mixing estimate, for functions on $X$ instead of on $T^1(X)$. Similar estimates for the latter are available \cite[Theorem 2]{Matheus}, and it would be natural to use those to express an analogous version of \eqref{cross_term_eq} but with $\phi \in L^2(T^1(X))$ given by 
    \begin{equation*}
        \phi(g) := \int_{0}^T K_{r, R}(\widetilde{\pi}\circ \widetilde{\G}_t(g), w) \dd t - T \frac{\mu(A_{r, R})}{\mu(\F)}.
    \end{equation*}
    The $L^2$-norm estimate is essentially the rest of the proof of \autoref{random_thm}, and one would be able to gain an extra factor of $(R-r) \log\left(\frac{2R}{R-r}\right)$ in the error term coming from cutting the geodesic into small pieces of length $T$. The issue is that \cite[Theorem 2]{Matheus} requires estimates for Lie derivatives of order $3$ in the angular direction, which would add a lot of complexity at diminishing returns. 
    
    Instead we are satisfied with the mild restriction on $R-r$, which is already enough to accommodate for instance $R-r \geq \exp(-R^{-1+\delta})$ for any fixed $\delta > 0$.
\end{remark}


\section{Variance for closed geodesics}

\subsection{Predictions from the random model}

We can rewrite \eqref{variance_eq} as
\begin{equation*}
    \Var(r, R; \Lambda_D) := \int_X \left( \sum_{\mathcal{C}\in \Lambda_D} \ell(\mathcal{C} \cap A_{r, R}(w)) - \frac{\mu(A_{r, R})}{\widetilde{\mu}(X)} \sum_{\mathcal{C}\in \Lambda_D} \ell(\mathcal{C}) \right)^2 \frac{d\widetilde{\mu}(w)}{\widetilde{\mu}(X)}.
\end{equation*}

If $J := (\sqrt{D}) \in \Cl_D^+$, then the closed geodesic corresponding to any $B \in \Cl_D^+$ is the same as that corresponding to $JB^{-1}$, but with opposite orientation. If $B^2 = J$, they are the same and correspond to a (so-called reciprocal \cite{Reciprocal}) closed geodesic that goes through its image twice, once in each orientation. Let $a_D := |\{B \in \Cl_D^+ : B^2 = J\}|$, so $h^+_D = a_D + 2 b_D$. The images of the closed geodesics from $\Lambda_D$ in $\Gamma \backslash \H$ correspond to $a_D$ geodesic segments of length $\log{\eps_D^+}$ with multiplicity $2$, and $b_D$ geodesic segments of length $2 \log{\eps_D^+}$ also with multiplicity $2$. Furthermore, the height of each of those closed geodesics in the fundamental domain is $\leq \sqrt{D}/2$ \cite[Proposition 3.1]{ELMV}. Therefore, if we model each of those geodesic segments using independent (except for the multiplicities) random geodesic segments in $X$, with a cutoff $A > \sqrt{D}/2$, we may expect

\begin{equation*}
    \begin{split}
        \Var(r, R; \Lambda_D) & \approx 4\left(a_D \Var_A(r, R; \log{\eps_D^+}) + b_D \Var_A(r, R; 2\log{\eps_D^+}) \right) \\
        & \sim 4\left(a_D \frac{16 \cdot \log{\eps_D^+} \cdot R^3}{\pi} \mathbf{G}\left(\frac{r}{R}\right) + b_D \frac{16 \cdot 2\log{\eps_D^+} \cdot R^3}{\pi} \mathbf{G}\left(\frac{r}{R}\right) \right) \\
        & = \frac{64 \sqrt{D} L(1, \chi_D) R^3}{\pi} \mathbf{G}\left(\frac{r}{R}\right),
    \end{split}
\end{equation*}
at least for $\log\left(\frac{\sqrt{D}}{2}\right) < \log{A} = o\left(R^{-1} \log^{-1}\left(\frac{1}{R}\right)\right)$ (taking \autoref{mixing_rmk} into account, but already from \autoref{random_thm} for thick annuli with $R-r \gg R$). It would suffice to restrict to $R \leq (\log{D})^{-1-\delta}$ for any fixed $\delta > 0$. Therefore, being a bit conservative this leads to the conjecture below.

\begin{conjecture}\label{main_conj}
Let $\delta > 0$ be given. If $0 \leq r < R \leq D^{-\delta}$, then as $D \to \infty$ through squarefree fundamental discriminants,
\begin{equation*}
    \Var(r, R; \Lambda_D) \sim \frac{64 \sqrt{D} L(1, \chi_D) R^3}{\pi} \mathbf{G}\left(\frac{r}{R}\right).
\end{equation*}
\end{conjecture}

Our main result confirms the conjecture for sufficiently small annuli that are not too thin, and in particular for small balls. Observe that the allowed range of radii intersects the regime where one would expect equidistribution, i.e. $\mu(A_{r, R}) \geq D^{-1 + \delta}$.

\begin{theorem}\label{main_thm}
    Let $\delta > 0$ be given. If $0 \leq r < R \leq D^{-\frac{5}{12} - \delta}$ and $R-r \gg R$, then as $D \to \infty$ through squarefree fundamental discriminants,
\begin{equation*}
    \Var(r, R; \Lambda_D) \sim \frac{64 \sqrt{D} L(1, \chi_D) R^3}{\pi} \mathbf{G}\left(\frac{r}{R}\right).
\end{equation*}
\end{theorem}

\begin{remark}
    The restriction $R-r \gg R$ in \autoref{main_thm} is mostly technical in nature, due to the fact that the behavior of the weight function $h_{r, R}(t)$ changes when $R-r \ll R$. We stick with it for simplicity, since it is enough to cover the most interesting case of balls ($r = 0$).
\end{remark}

\begin{remark}
    The proof of \autoref{main_thm} actually gives a power-saving error term of the form $O_\eps(D^\frac{1}{2}R^{3 + \eps})$ for any $\eps>0$ sufficiently small (depending on $\delta$). 
\end{remark}

\subsection{Spectral expansion and automorphic transformations}

Let $D>0$ be a squarefree fundamental discriminant, $\chi_D$ be the primitive quadratic character modulo $D$, and $\mathcal{B}_0(\Gamma)$ be an orthonormal basis of the space of Maa{\ss} cusp forms for the modular group $\Gamma$, which we may choose to consist of Hecke–Maa{\ss} cusp forms. 

Expressing the variance in terms of the automorphic kernel $K_{r, R}$, performing a spectral expansion, and using the work of Duke--Imamo\={g}lu--T\'{o}th \cite{DIT} to compute the resulting Weyl sums, we are left with $L$-functions.

\begin{lemma}[Spectral expansion of the variance {\cite[Lemma 2.20]{HR19}}] \label{spectral_expansion_lemma}
    We have 
    \begin{equation*}
        \begin{split}
            \Var(r, R; \Lambda_D) =\ & \frac{\sqrt{D}}{2 \widetilde{\mu}(X)} \sum_{f \in \mathcal{B}_0(\Gamma)} \frac{L\left(\frac{1}{2}, f\right) L\left(\frac{1}{2}, f \otimes \chi_D\right)}{L\left(1, \sym^2 f\right)} H(t_f) \left| h_{r, R}(t_f) \right|^2 \\
            & + \frac{\sqrt{D}}{4 \pi \widetilde{\mu}(X)} \int_{-\infty}^\infty \left| \frac{\zeta\left(\frac{1}{2} + it\right) L\left( \frac{1}{2} + it, \chi_D \right)}{\zeta(1 + 2it)} \right|^2 H(t) \left| h_{r, R}(t) \right|^2 \dd t,
        \end{split}
    \end{equation*}
    where
    \begin{equation}\label{gamma_factor_def_eq}
        H(t) := \frac{\Gamma\left(\frac{1}{4} + \frac{it}{2}\right)^2 \Gamma\left(\frac{1}{4} - \frac{it}{2} \right)^2}{\Gamma\left(\frac{1}{2} + it\right) \Gamma\left(\frac{1}{2} - it \right)} = \frac{4\pi}{|t|+1} + O\left(\frac{1}{(|t|+1)^2} \right).
    \end{equation}
\end{lemma}

Before the next lemma we need to establish some notation. Recall that the Mellin transform $\widehat{W}$ of a function $W : (0, \infty) \to \C$ is given by
\begin{equation*}
    \widehat{W}(s) := \int_0^\infty W(x) x^s \frac{dx}{x}
\end{equation*}
for $s \in \C$ for which the integral is absolutely convergent, and conversely the inverse Mellin transform $\widecheck{\mathcal{W}}$ of a holomorphic function $\mathcal{W} : \{ s \in \C : a < \Im(s) < b\} \to \C$ is given by
\begin{equation*}
    \widecheck{\mathcal{W}}(x) := \frac{1}{2\pi i} \int_{\sigma-i\infty}^{\sigma + i\infty} \mathcal{W}(s) x^{-s} \dd s
\end{equation*}
for $a < \sigma < b$ and $x \in (0, \infty)$ for which the integral converges absolutely.

\begin{lemma}[Automorphic transformations {\cite[Corollary 5.7]{HR19}}] \label{automorphic_lemma}
    Let $h(t)$ be an even holomorphic function in the strip $-2M < \Im(t) < 2M$ for some $M \geq 20$ with zeroes at $\pm \left(n - \frac{1}{2}\right)i$ for $n \in \{1, 2, \dots , 2M\}$ and satisfying $h(t) \ll (|t| + 1)^{-2M}$ in this region. Then the moment
    \begin{equation*}
        \sum_{f \in \mathcal{B}_0(\Gamma)} \frac{L\left(\frac{1}{2}, f\right) L\left(\frac{1}{2}, f \otimes \chi_D\right)}{L\left(1, \sym^2 f\right)} h(t_f) + \frac{1}{2 \pi} \int_{-\infty}^\infty \left| \frac{\zeta\left(\frac{1}{2} + it\right) L\left( \frac{1}{2} + it, \chi_D \right)}{\zeta(1 + 2it)} \right|^2 h(t) \dd t
    \end{equation*}
    is equal to the sum of the main term
    \begin{equation}\label{general_main_term_eq}
        2 L(1, \chi_D) \int_{-\infty}^\infty h(t) \dspec t
    \end{equation}
    and the shifted convolution
    \begin{equation}\label{general_shifted_conv_eq}
        \begin{split}
            \frac{2}{\sqrt{D}} \sum_\pm \sum_{D_1 D_2 = D} \sum_{\substack{m = 1 \\ m \not= \mp D_2}}^\infty &\chi_1(\sgn(m \pm D_2)) \lambda_{\chi_1, \chi_2}(m, 0) \lambda_{\chi_1, \chi_2}(|m \pm D_2|, 0) \\
            &\times \frac{1}{2\pi i} \int_{\sigma_1 - i\infty}^{\sigma_1 + i\infty} \widehat{\mathscr{K}^- h}(s) \widehat{\mathcal{J}_0^\pm}(1-s) \left(\frac{m}{D_2}\right)^\frac{s-1}{2} d s,
        \end{split}
    \end{equation}
    where $1 - M < \sigma_1 < -1$, 
    \begin{equation*}
        \begin{split}
            \left(\mathscr{K}^- h\right)(x) := \int_{-\infty}^{\infty} h(t) \mathcal{J}^-_t(x) \dspec t &, \quad \quad \quad \quad \dspec t := \frac{1}{2 \pi^2} t \tanh\left(\pi t\right) \dd t, \\
            \mathcal{J}^-_t(x) := 4 \cosh\left(\pi t\right) K_{2 it}(4 \pi x) &, \quad \quad \quad \quad \mathcal{J}^+_0(x) := -2 \pi Y_0(4 \pi x),
        \end{split}
    \end{equation*}
    the decomposition $\chi_D = \chi_1 \chi_2$ corresponds to $D = D_1 D_2$, and
    \begin{equation*}
        \lambda_{\chi_1, \chi_2}(m, 0) := \sum_{ab = m} \chi_1(a) \chi_2(b).
    \end{equation*}
\end{lemma}

Combining this with work of M. Young \cite{Young} leads to the following bound for moments of $L$-functions, which will be useful for bounding some of our error terms later on.

\begin{lemma}[Dyadic moment bound {\cite[Proposition 2.35 (1)]{HR19}}] \label{moment_bound_lemma}
    For $T \geq 1$, we have
    \begin{equation*}
        \begin{split}
            \sum_{\substack{f \in \mathcal{B}_0(\Gamma) \\ T \leq t_f \leq 2T}} \frac{L\left(\frac{1}{2}, f \right) L\left(\frac{1}{2}, f \otimes \chi_D \right)}{L(1, \sym^2 f)} + \frac{1}{2\pi} \int\limits_{T \leq |t| \leq 2T} &\left|\frac{\zeta\left( \frac{1}{2} + it \right) L\left( \frac{1}{2} + it, \chi_D\right)}{\zeta(1+ 2it)}\right|^2 dt \\
            & \ll_\eps 
            \begin{cases}
                D^{\frac{1}{3} + \eps} T^{2+\eps} & \text{for } T \ll D^\frac{1}{12}, \\
                D^{\frac{1}{2} + \eps} & \text{for } D^\frac{1}{12} \ll T \ll D^\frac{1}{4}, \\
                D^\eps T^{2+\eps} & \text{for } T \gg D^\frac{1}{4}.
            \end{cases}
        \end{split}
    \end{equation*}
\end{lemma}

\subsection{Choice of test function} To prove \autoref{main_thm}, we will start with the expression in \autoref{spectral_expansion_lemma} and approximate the weights $H(t) |h_{r, R}(t)|^2$ by a function $h(t)$ satisfying the conditions of \autoref{automorphic_lemma}. The error terms coming from switching from one set of weights to the other may be bounded using \autoref{moment_bound_lemma}, and the problem will be reduced to evaluating the main term \eqref{general_main_term_eq} and the error term \eqref{general_shifted_conv_eq}. We once again follow \cite{HR19}, adapting their construction to our context.

The conditions of \autoref{automorphic_lemma} require that $h(t)$ be even, extend holomorphically to $|\Im(t)| < 2M$, have zeros at $\pm\left(n-\frac{1}{2}\right)i$ for $n \in \{ 1, 2, \dots, 2M\}$, and satisfy $h(t) \ll (|t| + 1)^{-2M}$ for some integer $M \geq 20$. From now on, fix a sufficiently large constant $M \in \N$. 

First we localize $h(t)$ to the region $[-T_2, -T_1]\cup[T_1, T_2]$, where $T_1 := R^{-1+\alpha}$ and $T_2:= R^{-1-\alpha}$ for a sufficiently small fixed constant $\alpha > 0$. This is because the main contribution to $\Var(r, R; \Lambda_D)$ will come from this range when $R-r \gg R$. To achieve this localization, let
\begin{equation*}
    h_1(t) := e^{-\left(\frac{t}{T_2}\right)^{2M}} \left( 1 - e^{-\left(\frac{t}{T_1}\right)^{2M}}\right),
\end{equation*}
which is even and for $|\Im(t)| < 2M$ satisfies
\begin{equation}\label{h_1_bound_eq}
    h_1(t) = 
    \begin{cases}
        O\left(\left(\frac{|\Re(t)| + 1}{T_1}\right)^{2M}\right) & \text{for } |\Re(t)| \leq T_1, \\
        1 + O\left(\left(\frac{\Re(t)}{T_2}\right)^{2M} + e^{-\left(\frac{\Re(t)}{T_1}\right)^{2M}}\right) & \text{for } T_1 \leq |\Re(t)| \leq T_2,\\
        O\left(e^{-\left(\frac{\Re(t)}{T_2 }\right)^{2M}}\right) & \text{for } |\Re(t)| \geq T_2.
    \end{cases}
\end{equation}
Moreover, for $j \in \{1, 2, \dots, 2M\}$ and $t \in \R$,
\begin{equation}\label{h_1_deriv_bound_eq}
    h_1^{(j)}(t) \ll_j 
    \begin{cases}
        \frac{|t|^{2M-j}}{T_1^{2M}} & \text{for } |t| \leq T_1, \\
        \frac{|t|^{2M-j}}{T_2^{2M}} + \frac{|t|^{(2M - 1)j}}{T_1^{2Mj}} e^{-\left(\frac{t}{T_1}\right)^{2M}} & \text{for } T_1 \leq |t| \leq T_2, \\
        \frac{|t|^{(2M - 1)j}}{T_2^{2Mj}} e^{-\left(\frac{t}{T_2}\right)^{2M}}  & \text{for } |t| \geq T_2.
    \end{cases}
\end{equation}

Next, ignoring the Bessel factors for now, we see from \autoref{SH_asymp_lemma} and \eqref{gamma_factor_def_eq} that a factor asymptotic to $16 \pi^3/|t|^3$ arises. Therefore consider
\begin{equation*}
    h_2(t) := 2 (2 \pi)^{-4M + 1} (4M+3)^{-3} \frac{\Gamma\left(\frac{2M}{4M+3} + \frac{it}{4M+3} \right)^{4M+3} \Gamma\left(\frac{2M}{4M+3} - \frac{it}{4M+3} \right)^{4M+3}}{\Gamma\left(\frac{1}{2} + it \right) \Gamma\left(\frac{1}{2} - it \right)},
\end{equation*}
which is even and holomorphic in the strip $|\Im(t)| < 2M$, where it has zeros at $\pm\left(n-\frac{1}{2}\right)i$ for $n \in \{ 1, 2, \dots, 2M\}$ and satisfies
\begin{equation}\label{h_2_bound_eq}
    h_2(t) = \frac{16 \pi^3}{(|t| + 1)^3} + O\left(\frac{1}{(|t| + 1)^4}\right),
\end{equation}
by Stirling's formula. Furthermore, for $j \in \Z_{\geq 0}$ and $t \in \R$,
\begin{equation}\label{h_2_deriv_bound_eq}
    h_2^{(j)}(t) \ll_j (|t| + 1)^{-j-3}.
\end{equation}

Finally, let
\begin{equation*}
    h_3(t) := (R \cdot J_1(Rt) - r \cdot J_1(rt))^2,
\end{equation*}
which is entire (as this is the case for $J_1$) and even (as $J_1$ is odd \cite[8.476.1]{Integrals}). Using the crude bound
\begin{equation}\label{J_1_crude_bound_eq}
    J_1(z) \ll
    \begin{cases}
        |z| & \text{for } |z| \leq 1, \\
        \frac{e^{|\Im(z)|}}{\sqrt{|z|}} & \text{for } |z| \geq 1
    \end{cases}
\end{equation}
\cite[8.411.3 and 8.451.1]{Integrals}, for $|\Im(t)| < 2M$ we get
\begin{equation}\label{h_3_bound_eq}
    h_3(t) \ll 
    \begin{cases}
        R^4 |t|^2 & \text{for } |t| \leq \frac{1}{R},\\
        \frac{R}{|t|} & \text{for } |t| \geq \frac{1}{R}.
    \end{cases}
\end{equation}

Also, for $j \in \Z_{\geq 0}$ and $y \in \R$ we have
\begin{equation*}
    J_1^{(j)}(y) \ll_j
    \begin{cases}
        |y| & \text{for } |y|\leq 1 \text{ and } j \text{ even},\\
        1 & \text{for } |y|\leq 1 \text{ and } j \text{ odd},\\
        \frac{1}{\sqrt{|y|}} & \text{for } |y| \geq 1
    \end{cases}
\end{equation*}
\cite[8.471.2, 8.411.4, and 8.451.1]{Integrals}, which for $t \in \R$ gives
\begin{equation}\label{h_3_deriv_bound_eq}
    h_3^{(j)}(t) \ll_j
    \begin{cases}
        R^4 |t|^{2-j} & \text{for } |t| \leq \frac{1}{R} \text{ and } j \in \{0, 1\},\\
        R^{2+j} & \text{for } |t| \leq \frac{1}{R} \text{ and } j\geq 2,\\
        \frac{R^{1+j}}{|t|} & \text{for } |t| \geq \frac{1}{R}.
    \end{cases}
\end{equation}

We choose the test function 
\begin{equation}\label{h_def_eq}
    h(t) := h_1(t) h_2(t) h_3(t),
\end{equation}
so that combining \eqref{h_1_bound_eq}, \eqref{h_2_bound_eq}, and \eqref{h_3_bound_eq} gives the following upper bounds and asymptotics for $h$.

\begin{lemma}\label{h_bound_lemma}
    For $|\Im(t) | < 2M$,
    \begin{equation}\label{h_bound_eq}
        h(t) \ll
        \begin{cases}
            \frac{R^4 (|\Re(t)| + 1)^{2M-1}}{T_1^{2M}} &\text{for } |t| \leq T_1,\\
            \frac{R^4}{|\Re(t)|} & \text{for } T_1 \leq |t| \leq \frac{1}{R},\\
            \frac{R}{|\Re(t)|^4} & \text{for } \frac{1}{R} \leq |t| \leq T_2,\\
            \frac{R}{|\Re(t)|^4} e^{-\left(\frac{\Re(t)}{T_2}\right)^{2M}} & \text{for } |t| \geq T_2.
        \end{cases}
    \end{equation}
    Furthermore, if $t \in \R$ then 
    \begin{equation}
        h(t) = \frac{4 \pi}{|t|} \left( 2 \pi \frac{R \cdot J_1(Rt) - r\cdot J_1(rt)}{t} \right)^2 + 
        \begin{cases}
            O\left(\frac{R^4}{|t|^2} + \frac{R^4 |t|^{2M-1}}{T_2^{2M}} + \frac{R^4}{|t|}e^{-\left(\frac{t}{T_1}\right)^{2M}}\right) & \text{for } T_1 \leq |t| \leq \frac{1}{R},\\
            O\left(\frac{R}{|t|^5} + \frac{R |t|^{2M-4}}{T_2^{2M}} \right) & \text{for } \frac{1}{R} \leq |t| \leq T_2. 
        \end{cases}
    \end{equation}
\end{lemma}

We record the following important definitions and bounds for future reference:
\begin{equation}\label{reference_eq}
    R \ll D^{-\frac{5}{12} - \delta}, \quad \quad R-r \gg R, \quad \quad T_1 = R^{-1+\alpha}, \quad \quad T_2 = R^{-1 - \alpha},
\end{equation}
where $\alpha, \delta > 0$ are sufficiently small fixed constants and $M \in \N$ is a sufficiently large fixed constant.

\subsection{Change of test function for the variance}
    
\begin{lemma}
    Under the assumptions of \eqref{reference_eq} and for $h(t)$ as in \eqref{h_def_eq}, we have 
    \begin{equation}\label{change_of_var_eq}
        \begin{split}
            \Var(r, R; \Lambda_D) =\ & \frac{\sqrt{D}}{2 \widetilde{\mu}(X)} \sum_{f \in \mathcal{B}_0(\Gamma)} \frac{L\left(\frac{1}{2}, f\right) L\left(\frac{1}{2}, f \otimes \chi_D\right)}{L\left(1, \sym^2 f\right)} h(t_f) \\
            & + \frac{\sqrt{D}}{4 \pi \widetilde{\mu}(X)} \int_{-\infty}^\infty \left| \frac{\zeta\left(\frac{1}{2} + it\right) L\left( \frac{1}{2} + it, \chi_D \right)}{\zeta(1 + 2it)} \right|^2 h(t) \dd t \\
            & + O_\eps\left(D^{\frac{1}{2}+\eps} R^{3-\eps} \left( R D^{\frac{5}{12}} + R T_1 + \frac{R D^{\frac{1}{2}}}{T_1} + \frac{1}{(R T_2)^2} \right) \right).
        \end{split}
    \end{equation}
\end{lemma}

\begin{proof}
    Follows from the spectral expansion in \autoref{spectral_expansion_lemma} and a change of test function. The error term is estimated using the bounds and asymptotics in \autoref{SH_asymp_lemma}, \eqref{gamma_factor_def_eq}, and \autoref{h_bound_lemma}, considering each of the ranges separately. More specifically, if we denote $\Delta(t) := H(t) |h_{r, R}(t)|^2 - h(t)$, then putting those bounds together yields, for $t \in \R$,
    \begin{equation*}
        \Delta(t) \ll 
        \begin{cases}
            \frac{R^4}{1+|t|} & \text{for } |t| \leq T_1,\\
            \frac{R^4}{|t|^2} + \frac{R^4|t|^{2M-1}}{T_2^{2M}} + \frac{R^4}{|t|}e^{-\left(\frac{t}{T_1}\right)^{2M}} & \text{for } T_1 \leq |t| \leq \frac{1}{R},\\
            \frac{R}{|t|^5} +\frac{R|t|^{2M-4}}{T_2^{2M}} & \text{for } \frac{1}{R}\leq |t|\leq T_2,\\
            \frac{R}{|t|^4} & \text{for } |t|\geq T_2.
        \end{cases}
    \end{equation*}
    We used the fact that $\alpha > 0$ is small in the inequality above. 
    
    Therefore, combining this with \autoref{moment_bound_lemma} gives, for $T \geq 1$,
    \begin{equation*}
        \begin{split}
            \sum_{\substack{f \in \mathcal{B}_0(\Gamma) \\ T \leq t_f \leq 2T}} \frac{L\left(\frac{1}{2}, f\right) L\left(\frac{1}{2}, f \otimes \chi_D\right)}{L\left(1, \sym^2 f\right)} |\Delta(t_f)| + \frac{1}{2 \pi} \int\limits_{T \leq |t| \leq 2T} \left| \frac{\zeta\left(\frac{1}{2} + it\right) L\left( \frac{1}{2} + it, \chi_D \right)}{\zeta(1 + 2it)} \right|^2 |\Delta(t)| \dd t & \\
            \ll_\eps
            \begin{cases}
                D^{\frac{1}{3}+\eps} T^{1+\eps} R^4 & \text{for } 1 \leq T \ll D^\frac{1}{12},\\
                D^{\frac{1}{2}+\eps} \frac{R^4}{T} & \text{for } D^\frac{1}{12} \ll T \ll D^\frac{1}{4},\\
                D^\eps T^{1+\eps} R^4 & \text{for } D^\frac{1}{4} \ll T \leq T_1,\\
                D^\eps T^\eps \left( R^4 + R^4 T\left(\frac{T}{T_2}\right)^{2M} + R^4T e^{-\left(\frac{T}{T_1}\right)^{2M}} \right) & \text{for } T_1\leq T \leq \frac{1}{R},\\
                D^\eps T^\eps \left( \frac{R}{T^3} + \frac{R}{T^2} \left(\frac{T}{T_2}\right)^{2M} \right) & \text{for } \frac{1}{R} \leq T \leq T_2,\\
                D^\eps T^\eps \frac{R}{T^2} & \text{for }  T \geq T_2.
            \end{cases}&
        \end{split}
    \end{equation*}
    Here we recall that $T_1 = R^{-1+\alpha} \gg D^{(1-\alpha)\left(\frac{5}{12} + \delta\right)} \geq D^\frac{1}{4}$ for sufficiently small $\alpha>0$. Multiplying by $\sqrt{D}$ and summing over $T=2^k$ for $k\in \Z_{\geq 0}$ gives the claimed error term. 
    
\end{proof}

Observe that by \eqref{reference_eq} the error term is $O_\eps(D^{\frac{1}{2}} R^{3+\eps})$ for $\eps > 0$ sufficiently small, and therefore it is asymptotically smaller than the main term of \autoref{main_thm}, as $L(1, \chi_D) \gg_\eps D^{-\eps}$. 

\begin{remark}
    The error term $O_\eps(D^{\frac{11}{12}+\eps} R^{4-\eps})$ in \eqref{change_of_var_eq} is the only point in the proof of \autoref{main_thm} where the range $R \leq D^{-\frac{5}{12}-\delta}$ is tight. Instead of using the bound $\ll_\eps D^{\frac{1}{3}+\eps}T^{2+\eps}$ (coming from Young's work \cite{Young}) for the range $T\ll D^\frac{1}{12}$ of \autoref{moment_bound_lemma}, we could have tried to use the weaker bound $\ll_\eps D^{\frac{1}{2}+\eps}$ (which holds in this range by the argument in Humphries-Radziwi{\l}{\l} \cite[Proposition 2.35]{HR19}). This would produce a corresponding error term of size $O_\eps(D^{1+\eps}R^4)$ in \eqref{change_of_var_eq}, which is enough to obtain asymptotics for the variance if $R \leq D^{-\frac{1}{2} - \delta}$. Here it becomes clear that it is precisely the range of (conjectured) equidistribution, i.e. $R\geq D^{-\frac{1}{2}+\delta}$, which requires deeper arithmetic inputs.
\end{remark}

Applying \autoref{automorphic_lemma} to the first two terms of \eqref{change_of_var_eq}, we obtain the main term
\begin{equation}\label{main_term_eq}
    \frac{\sqrt{D}L(1, \chi_D)}{\widetilde{\mu}(X)} \int_{-\infty}^\infty h(t) \dspec t,
\end{equation}
where $\dspec t := \frac{1}{2 \pi^2} t \tanh\left(\pi t\right) \dd t$ as before, plus the shifted convolution
\begin{equation}\label{shifted_conv_eq}
    \begin{split}
        \frac{1}{\widetilde{\mu}(X)} \sum_\pm \sum_{D_1 D_2 = D} \sum_{\substack{m = 1 \\ m \not= \mp D_2}}^\infty &\chi_1(\sgn(m \pm D_2)) \lambda_{\chi_1, \chi_2}(m, 0) \lambda_{\chi_1, \chi_2}(|m \pm D_2|, 0) \\
        &\times \frac{1}{2\pi i} \int_{\sigma_1 - i\infty}^{\sigma_1 + i\infty} \widehat{\mathscr{K}^- h}(s) \widehat{\mathcal{J}_0^\pm}(1-s) \left(\frac{m}{D_2}\right)^\frac{s-1}{2} ds.
    \end{split}
\end{equation}

\subsection{Asymptotics for main term}

\begin{lemma}\label{main_term_lemma}
    Under the assumptions of \eqref{reference_eq} and for $h(t)$ as in \eqref{h_def_eq}, the main term \eqref{main_term_eq} is equal to
    \begin{equation*}
        \frac{64 \sqrt{D} L(1, \chi_D) R^3}{\pi} \mathbf{G}\left(\frac{r}{R}\right) + O_\eps\left(D^{\frac{1}{2} + \eps} \left( R^4 T_1 + \frac{R}{T_2^2} \right) \right).
    \end{equation*}
\end{lemma}

\begin{proof}
    Using the bounds and asymptotics of \autoref{h_bound_lemma}, combined with the fact that $h$ is even and the bound $L(1, \chi_D) \ll \log{D}$, we see that \eqref{main_term_eq} is equal to
    \begin{equation}\label{main_term_integral_range_eq}
        \frac{16 \pi \sqrt{D}L(1, \chi_D)}{\widetilde{\mu}(X)} \int_{T_1}^{T_2} \left( \frac{R\cdot J_1(Rt) - r \cdot J_1(rt)}{t} \right)^2 dt + O_\eps\left(D^{\frac{1}{2} + \eps} \left( R^4 T_1 + \frac{R}{T_2^2} \right) \right).
    \end{equation}
    Then \eqref{J_1_crude_bound_eq} allows us to complete the integral to $(0, \infty)$ under the same error term as above. 
    
    From \cite[6.574.2]{Integrals} it follows that
    \begin{equation}\label{J_square_term_eq}
        R^2 \int_0^\infty \frac{J_1(Rt)^2}{t^2} \dd t = \frac{R^3 \cdot \Gamma\left(\frac{1}{2}\right)}{4 \cdot \Gamma\left(\frac{3}{2}\right) \Gamma\left(\frac{5}{2}\right) \Gamma\left(\frac{3}{2}\right)} = \frac{4 R^3}{3 \pi},
    \end{equation}
    and similarly for the term corresponding to $r$. The cross-term can be evaluated using \cite[6.574.3]{Integrals}, which gives
    \begin{equation*}
        2 R r \int_0^\infty \frac{J_1(Rt)J_1(rt)}{t^2} \dd t = R r^2 \cdot {}_2F_1\left( \frac{1}{2}, -\frac{1}{2}; 2; \frac{r^2}{R^2} \right),
    \end{equation*}
    where ${}_2F_1$ denotes the ordinary hypergeometric function. By \cite[8.113.1, 8.114.1, and 9.137.14]{Integrals} we deduce that
    \begin{equation*}
        {}_2F_1\left( \frac{1}{2}, -\frac{1}{2}; 2; z^2 \right) = \frac{4}{3 \pi z^2} \left((1+z^2) \mathbf{E}(z) - (1-z^2) \mathbf{K}(z) \right).
    \end{equation*}
    Therefore,
    \begin{equation}\label{J_cross_term_eq}
        2 R r \int_0^\infty \frac{J_1(Rt)J_1(rt)}{t^2} \dd t = \frac{4 R^3}{3 \pi} \left(\left(1+\frac{r^2}{R^2}\right) \mathbf{E}\left(\frac{r}{R}\right) - \left(1-\frac{r^2}{R^2}\right) \mathbf{K}\left(\frac{r}{R}\right) \right).
    \end{equation}
    
    Combining \eqref{J_square_term_eq} and \eqref{J_cross_term_eq}, we conclude that
    \begin{equation*}
        \int_0^\infty \left( \frac{R \cdot J_1(Rt) - r \cdot J_1(rt)}{t} \right)^2 dt = \frac{4 R^3}{3 \pi} \mathbf{G} \left(\frac{r}{R}\right),
    \end{equation*}
    which gives the desired result.
    
\end{proof}

By \eqref{reference_eq}, the error term in \autoref{main_term_lemma} is $O_\eps(D^{\frac{1}{2}} R^{3+\eps})$ for $\eps > 0$ sufficiently small, so it is once again asymptotically smaller than the main term of \autoref{main_thm}.

\subsection{Bounds for shifted convolution}

To finish the proof of \autoref{main_thm}, it suffices to show that the shifted convolution \eqref{shifted_conv_eq} is asymptotically smaller than the main term obtained in the previous subsection. This requires considerably more work and involves a more careful consideration of the oscillatory behavior of the test function $h(t)$. The final result is indicated in the lemma below.

\begin{lemma}\label{shifted_conv_lemma}
    Under the assumptions of \eqref{reference_eq} and for $h(t)$ as in \eqref{h_def_eq}, the shifted convolution \eqref{shifted_conv_eq} is $O_\eps(D^{\frac{1}{2}} R^{3+\eps})$ for every $\eps > 0$ sufficiently small.
\end{lemma}

\begin{proof}
    We once again follow \cite{HR19}, with necessary modifications due to the fact that $D > 0$ and also the presence of oscillations coming from a Bessel function, instead of a trigonometric function, in our choice of $h(t)$. 
    
    By Mellin inversion -- where we use the convolution identity \cite[(A.6)]{Ivi} -- and the divisor bound, it suffices to show that
    \begin{equation}\label{error_expression_eq}
        \sum_\pm \sum_{D_2 | D} \sum_{m=1}^\infty m^\eps \left| \int_0^\infty (\mathscr{K}^- h) (x) \mathcal{J}_0^\pm \left( \sqrt{\frac{m}{D_2}}x\right)\dd x \right|
    \end{equation}
    is $O_\eps(D^{\frac{1}{2}} R^{3+\eps})$ for every $\eps > 0$ sufficiently small. We consider two different ranges for $m$.

    \begin{itemize}[leftmargin=*, align=left]
        \item[\textit{Case 1:}] $m > \sqrt{D_2}$.
        
        Via integration by parts and \cite[8.472.1 and 8.486.12]{Integrals}, the integral in \eqref{error_expression_eq} is
        \begin{equation}\label{error_integral_eq}
            \frac{c^\pm D_2}{16 \pi^2 m}\int_0^\infty \frac{1}{x^2} \mathscr{L}(x) B^\pm_2\left(4\pi\sqrt{\frac{m}{D_2}}x\right)\dd x,
        \end{equation}
        where 
        \begin{equation*}
            c^+ := -2\pi, \quad \quad c^- := 4, \quad \quad B^+_k(x) := Y_k(x), \quad \quad B^-_k(x) := K_k(x),
        \end{equation*}
        and
        \begin{equation*}
            \mathscr{L}(x) := 3(\mathscr{K}^-h)(x) -3x (\mathscr{K}^-h)'(x) + x^2(\mathscr{K}^-h)''(x).
        \end{equation*}
        We will split the integral in \eqref{error_integral_eq} into three different ranges for $x$ and bound each one separately.
        
        \begin{itemize}[leftmargin=*, align=left]
            \item[\textit{Sub-case 1a:}] $0 < x \leq 1$.
            
            By \cite[(A.2) and (A.4)]{BLM},
            \begin{equation*}
                \frac{d^j}{dx^j}\mathcal{J}^-_t(x) = \frac{(2\pi)^j \pi i}{\sinh(\pi t)} \sum_{n=0}^j \binom{j}{n} (I_{2it - j +2n}(4 \pi x) - I_{-2it - j + 2n}(4 \pi x)).
            \end{equation*}
            Combining this with the bound
            \begin{equation*}
                e^{-\pi |t|} I_{2it - j + 2n}(4 \pi x) \ll_{\Im(t), j} \frac{x^{-j+2(n - \Im(t))}}{(|\Re(t)|+1)^{\frac{1}{2}-j+2(n-\Im(t))}}
            \end{equation*}
            valid for $0 < x \ll \sqrt{|t|+1}$, which follows from a slight adaptation of \cite[(A.6)]{BLM}, we can shift contours to obtain
            \begin{equation*}
                x^j \frac{d^j}{dx^j} (\mathscr{K}^-h)(x) \ll_j \sum_{\pm} \sum_{n=0}^j x^{2(n-c_n)} \int\limits_{\Im(t) = \pm c_n} |h(t)| (|\Re(t)| + 1)^{j -2(n-c_n) + \frac{1}{2}} \dd t
            \end{equation*}
            for any choice of integers $-2M < c_n < 2M$ (observe that the poles of $\cosh^{-1}(\pi t)$ are cancelled by the zeros of $h(t)$). Choose $c_n = n - 2M + 1$ and apply \eqref{h_bound_eq} to conclude that for $0 < x \leq 1$, 
            \begin{equation}\label{L_bound_1}
                \mathscr{L}(x) \ll \frac{R^4 x^{4M-2}}{T_1^{2M}}.
            \end{equation}
            
            For future reference, we note that if $k \in \Z_{\geq 0}$ and $x \in \R_{> 0}$ then one has the general bound
            \begin{equation}\label{general_bessel_bound_eq}
                B_k^\pm(x) \ll_{k, \eps} 
                \begin{cases}
                    x^{-k-\eps} & \text{for } 0< x \leq 1, \\
                    \frac{1}{\sqrt{x}} & \text{for } x \geq 1
                \end{cases}
            \end{equation}
            for both Bessel functions in question \cite[Proposition 9]{HM}. Therefore, we conclude that the contribution of $0< x\leq 1$ to \eqref{error_integral_eq} is 
            \begin{equation*}
                \ll_\eps \frac{R^4}{T_1^{2M}}\left( \left(\frac{D_2}{m}\right)^{2+\eps} +  \left(\frac{D_2}{m}\right)^\frac{5}{4}\right).
            \end{equation*}
            Summing over $m > \sqrt{D_2}$ and $D_2|D$, this sub-case contributes $O_\eps(R^4 T_1^{-2M} D^{\frac{3}{2}+\eps})$ to \eqref{error_expression_eq}, which is $O_\eps(D^{\frac{1}{2}} R^{3+\eps})$ by \eqref{reference_eq}.
            
            \item[\textit{Sub-case 1b:}] $x \geq T_2 \log{T_2}$.
            
            We use \cite[(A.1)]{BLM} to write
            \begin{equation*}
                \frac{d^j}{dx^j}\mathcal{J}^-_t(x) = (-2\pi)^j \sum_{n=0}^j \binom{j}{n} 4 \cosh(\pi t) K_{2it - j + 2n}(4 \pi x)
            \end{equation*}
            and apply the uniform bound
            \begin{equation*}
                \cosh(\pi t) K_{2it -j + 2n}(4\pi x) \ll_{\Im(t), j} e^{\min\{0, -\pi(4x - |\Re(t)|)\}} \left( \frac{1+|\Re(t)| + 4\pi x}{4 \pi x}\right)^{|2\Im(t) + j-2n| + \frac{1}{10}}
            \end{equation*}
            valid for all $t\in \C$ \cite[(A.3)]{BLM}. Combining this with \eqref{h_bound_eq} gives, for $t \in \R$ and $x \geq T_2$,
            \begin{equation*}
                h(t) t  \frac{d^j}{dx^j}\mathcal{J}^-_t(x) \ll_j
                \begin{cases}
                    \frac{R^4 (|t|+1)^{2M}}{T_1^{2M}}e^{\pi|t|} e^{-4\pi x} & \text{for } |t|\leq T_1,\\
                    R^4 e^{\pi|t|} e^{-4\pi x} & \text{for } T_1\leq |t|\leq \frac{1}{R},\\
                    \frac{R}{|t|^3} e^{\pi|t|} e^{-4\pi x} & \text{for } \frac{1}{R}\leq |t|\leq T_2,\\
                    \frac{R}{|t|^3}e^{-\left(\frac{t}{T_2}\right)^{2M}} e^{\pi|t|} e^{-4\pi x} & \text{for } T_2\leq |t|\leq 4x,\\
                    \frac{R}{|t|^3}e^{-\left(\frac{t}{T_2}\right)^{2M}} \left(\frac{|t|}{x}\right)^{j+\frac{1}{10}} & \text{for }  |t|\geq 4x.
                \end{cases}
            \end{equation*}
            Considering each range separately (and in fact dividing the fourth range into $|t|\leq 2x$ and $|t|\geq 2x$) leads to
            \begin{equation}\label{L_bound_2}
                \begin{split}
                    \mathscr{L}(x) & \ll \sum_{j=0}^2 x^j \int_{-\infty}^\infty \left| h(t) t \frac{d^j}{dx^j}\mathcal{J}^-_t(x) \right| \dd t \\
                    &\ll \sum_{j=0}^2 x^j \left( R^4 e^{\pi T_2} e^{-4\pi x} + \frac{R}{T_2^2} e^{-2\pi x} + \frac{R}{x^3}e^{-\left(\frac{2x}{T_2}\right)^{2M}} + \frac{R T_2}{x^3}e^{-\left(\frac{4x}{T_2}\right)^{2M}} \right) \\
                    &\ll R^3 x^2 e^{-2\pi x} + \frac{R T_2}{x}e^{-\left(\frac{2x}{T_2}\right)^{2M}}.
                \end{split}
            \end{equation}
            
            Therefore, using \eqref{general_bessel_bound_eq} once again, the contribution of $x \geq T_2 \log{T_2}$ to the integral \eqref{error_integral_eq} is $O_A(T_2^{-A} (D_2/m)^\frac{5}{4})$ for any $A > 0$, which easily gives the desired bound of $O_\eps(D^{\frac{1}{2}} R^{3+\eps})$ for the corresponding contribution to \eqref{error_expression_eq}.
            
            \item[\textit{Sub-case 1c:}] $1 < x < T_2 \log{T_2}$.
            
            We use the identity
            \begin{equation*}
                (\mathscr{K}^-h)(x) = \frac{1}{\pi} \int_{-\infty}^\infty e(2x \sinh(\pi u)) \int_{-\infty}^\infty h(t) t \tanh(\pi t) e(-ut) \dd  t \dd u
            \end{equation*}
            from \cite[(A.8)]{BLM}, which is valid due to the rapid decay of $h(t)$, following from \eqref{h_bound_lemma}. Then integrating by parts in $u$ gives
            \begin{equation}\label{L_fourier_eq}
                \mathscr{L}(x) = \frac{1}{\pi} \int_{-\infty}^\infty e(2 x \sinh(\pi u)) \int_{-\infty}^\infty \widetilde{h}(t) (c_0(u) + c_1(u) t + c_2(u) t^2) e(-ut)\dd t \dd u,
            \end{equation}
            where
            \begin{equation*}
                \widetilde{h}(t) := h(t) t \tanh(\pi t)
            \end{equation*}
            and
            \begin{equation*}
                \begin{split}
                    c_0(u) &:= 8 - 8\tanh^2(\pi u) + 3 \tanh^4(\pi u), \\
                    c_1(u) &:= -14 i \tanh(\pi u) + 6i \tanh^3(\pi u),\\
                    c_2(u) &:= -4 \tanh^2(\pi u).
                \end{split}
            \end{equation*}
            
            From $\frac{d}{dt}\tanh(\pi t) = \pi \sech(\pi t)^2$ and $\frac{d}{dt}\sech(\pi t) = -\pi \tanh(\pi t) \sech(\pi t)$ we can show by induction that for $j \geq 1$, there is a polynomial $Q_j(x, y)$ such that
            \begin{equation*}
                \frac{d^j}{dt^j}\tanh(\pi t) = \sech(\pi t)^2 \cdot Q_j(\tanh(\pi t), \sech(\pi t)) \ll_j e^{-2\pi |t|}, 
            \end{equation*}
            which will be negligible in what follows. Combining such a bound with \eqref{h_1_deriv_bound_eq}, \eqref{h_2_deriv_bound_eq}, and \eqref{h_3_deriv_bound_eq} we conclude that for $j \in \{0, 1, \dots, 2M\}$ and $t \in \R$,
            \begin{equation}\label{h_tilde_der_bound_eq}
                \begin{split}
                    \widetilde{h}^{(j)}(t) & \ll_j
                    \begin{cases}
                        \left(\frac{t}{T_1}\right)^{2M} \cdot\frac{1}{(|t|+1)^3} \cdot R^4|t|^2 \cdot |t|\cdot \frac{1}{|t|^j} & \text{for } |t| \leq T_1,\\
                        \left(1+\left(\frac{t}{T_2}\right)^{2M} + \left(\frac{t}{T_1}\right)^{2Mj}e^{-\left( \frac{t}{T_1}\right)^{2M}} \right) \cdot \frac{1}{|t|^3}\cdot R^4 |t|^2 \cdot|t| \cdot \frac{1}{|t|^j} & \text{for } T_1 \leq |t| \leq \frac{1}{R},\\
                        \left(1+\left(\frac{t}{T_2}\right)^{2M} + \left(\frac{t}{T_1}\right)^{2Mj}e^{-\left( \frac{t}{T_1}\right)^{2M}} \right) \cdot \frac{1}{|t|^3}\cdot \frac{R}{|t|} \cdot|t| \cdot R^j & \text{for } \frac{1}{R} \leq |t| \leq T_2,\\
                        \left(\frac{t}{T_2}\right)^{2Mj} e^{-\left( \frac{t}{T_2} \right)^{2M}} \cdot \frac{1}{|t|^3}\cdot \frac{R}{|t|}\cdot |t|\cdot R^j  & \text{for } |t| \geq T_2
                    \end{cases}\\
                    &\ll_j
                    \begin{cases}
                        \frac{R^4 (|t| + 1)^{2M-j}}{T_1^{2M}} & \text{for } |t| \leq T_1,\\
                        \frac{R^4}{T_1^j} & \text{for } T_1 \leq |t| \leq \frac{1}{R},\\
                        \frac{R}{|t|^3 T_1^j} & \text{for } \frac{1}{R} \leq |t| \leq T_2,\\
                        \frac{R^{1+j}}{|t|^3}e^{-\frac{1}{2}\left(\frac{t}{T_2}\right)^{2M}} & \text{for } |t| \geq T_2.
                    \end{cases}
                \end{split}
            \end{equation}
            
            We now bound \eqref{L_fourier_eq} by dividing the integral over $u$ into the ranges $|u| \leq v$ and $|u| > v$, where $v \in (0, 1)$ will be chosen later. In the case $|u| > v$, we estimate the integral over $t$ by integrating by parts $2M$ times. Since $c_i(u) \ll |u|^i$ for $i \in \{0, 1, 2\}$, \eqref{h_tilde_der_bound_eq} shows that the contribution of this range to $\mathscr{L}(x)$ is
            \begin{equation*}
                \ll R \log(1/R) T_1^{-2M} v^{-2M+1}(R+v)^2.
            \end{equation*}
            
            For $|u| \leq v <1$, we Taylor expand twice to get $e(2x \sinh(\pi u)) = e(2x (\pi u + O(u^3))) = e(2 \pi xu) + O(xu^3)$, as long as $xv^3 < 1$ (which will be the case for our choice of $v$). Plugging this into \eqref{L_fourier_eq} and using \eqref{h_tilde_der_bound_eq}, the error term is
            \begin{equation*}
                \ll R \log(1/R) x v^4 (R+v)^2.
            \end{equation*}
            
            To make the two error terms collected so far match, we choose 
            \begin{equation*}
                v = T_1^{-1+\frac{3}{2M+3}}x^{-\frac{1}{2M+3}},
            \end{equation*}
            which satisfies the necessary restrictions since in the present sub-case $x < T_2 \log{T_2}$. In the remaining integral over $|u| \leq v$ we use
            \begin{equation*}
                \begin{split}
                    c_0(u) &= 8 + O(u^2),\\
                    c_1(u) &= -14 i \pi u + O(u^3),\\
                    c_2(u) &= -4\pi^2 u^2 + O(u^4).
                \end{split}
            \end{equation*}
            
            The contribution of these error terms to \eqref{L_fourier_eq} is 
            \begin{equation*}
                \ll R\log(1/R) v^3 (R+v)^2 \ll R \log(1/R) v^5,
            \end{equation*}
            as $1 < x < T_2 \log{T_2}$. Finally, we can complete the integral over $u$ to $(-\infty, \infty)$ under an error term 
            \begin{equation*}
                \ll R \log(1/R) T_1^{-2M} v^{-2M+1}(R+v)^2 = R \log(1/R) x v^4(R+v)^2 \ll R \log(1/R) x v^6,
            \end{equation*}
            by the argument via integration by parts from before. Therefore,
            \begin{equation}\label{L_fourier_bound_interm_eq}
                \begin{split}
                    \mathscr{L}(x) &= \frac{1}{\pi} \int_{-\infty}^\infty e(2\pi x u ) \int_{-\infty}^\infty \widetilde{h}(t) (8 -14i\pi u t -4 \pi^2 u^2 t^2) e(-ut) \dd t \dd u + O\left(R \log(1/R) v^5 (1+xv)\right)\\
                    & = \frac{1}{\pi}\left( 3 \widetilde{h}(2\pi x) - 3 x \frac{d}{dx}\left[ \widetilde{h}(2\pi x) \right] + x^2 \frac{d^2}{dx^2}\left[ \widetilde{h}(2\pi x) \right]\right) + O\left(R \log(1/R) v^5 (1+xv)\right),
                \end{split}
            \end{equation}
            where the double integral was evaluated via Fourier inversion. For $1 < x < T_2 \log{T_2}$ we have
            \begin{equation*}
                R \log(1/R) v^5 (1+xv) = R \log(1/R) \left( \frac{x^{-\frac{5}{2M+3}}}{T_1^{5 - \frac{15}{2M+3}}} + \frac{x^{1-\frac{5}{2M+3}}}{T_1^{6 - \frac{18}{2M+3}}}\right) \ll \frac{R}{T_1^\frac{9}{2}}
            \end{equation*}
            due to \eqref{reference_eq}. Applying this combined with \eqref{h_tilde_der_bound_eq} to \eqref{L_fourier_bound_interm_eq}, we obtain
            \begin{equation}\label{L_bound_3}
                \mathscr{L}(x) \ll \frac{R }{T_1^{\frac{9}{2}}} +
                \begin{cases}
                    \frac{R^4 x^2}{T_1^2} & \text{for } 1\leq x \leq \frac{1}{R},\\
                    \frac{R}{x T_1^2} & \text{for } \frac{1}{R} \leq x \leq T_2 \log{T_2}.
                \end{cases}
            \end{equation}
            Using the bound above and \eqref{general_bessel_bound_eq}, the contribution of $1 < x < T_2 \log{T_2}$ to \eqref{error_integral_eq} is
            \begin{equation*}
                \begin{split}
                    &\ll_\eps \frac{D_2}{m} \int_1^{T_2\log{T_2}} \frac{1}{x^2} |\mathscr{L}(x)| \left( \left(\frac{D_2}{m x^2}\right)^{1+\eps} + \left(\frac{D_2}{m}\right)^\frac{1}{4}\frac{1}{\sqrt{x}} \right) \dd x\\
                    &\ll \frac{R}{T_1^\frac{9}{2}} \left( \left(\frac{D_2}{m}\right)^{2+\eps} + \left(\frac{D_2}{m}\right)^\frac{5}{4} \right).
                \end{split}
            \end{equation*}
            Summing over $M > \sqrt{D_2}$ and $D_2|D$, this sub-case adds $O_\eps(RT_1^{-\frac{9}{2}}D^{\frac{3}{2}+\eps})$ to \eqref{error_expression_eq}. This is the most delicate range, but from \eqref{reference_eq} we see that it contributes $O_\eps(D^{\frac{1}{2}} R^{3+\eps})$, as desired.
        \end{itemize}
        
        \item[\textit{Case 2:}] $1 \leq m \leq \sqrt{D_2}$.
        
        In this case we directly bound the integral from \eqref{error_expression_eq}, which is
        \begin{equation}\label{error_integral_eq_2}
            c^\pm \int_0^\infty (\mathscr{K}^-h)(x) B_0^\pm\left(4 \pi \sqrt{\frac{m}{D_2}}x\right)\dd x.
        \end{equation}
         The strategy is to divide it into the same three ranges for $x$, and observe that the bounds \eqref{L_bound_1}, \eqref{L_bound_2}, and \eqref{L_bound_3} for $\mathscr{L}(x)$ are actually bounds for $\max_{j \in \{0, 1, 2\}} \left|x^j \frac{d^j}{dx^j} (\mathscr{K}^-h)(x)\right|$, so they apply verbatim to $(\mathscr{K}^-h)(x)$. We simply combine them with \eqref{general_bessel_bound_eq} to estimate \eqref{error_integral_eq_2}.
         
        \begin{itemize}[leftmargin=*, align=left]
            \item[\textit{Sub-case 2a:}] $0 < x \leq 1$.
            
            We see from \eqref{L_bound_1} that the contribution of $0 < x \leq 1$ to \eqref{error_integral_eq_2} is $O(R^4 T_1^{-2M} (D_2/m)^\frac{1}{4})$,            so this corresponds to a term of size $O_\eps(R^4 T_1^{-2M}D^{\frac{5}{8}+\eps})$ in \eqref{error_expression_eq}, which is acceptable.
            
            \item[\textit{Sub-case 2b:}] $x \geq T_2 \log{T_2}$.
            
            From \eqref{L_bound_2}, the contribution of this sub-case to \eqref{error_integral_eq_2} is $O_A(T_2^{-A} (D_2/m)^\frac{1}{4})$ for any $A > 0$, and this easily leads to an acceptable error term of $O_{A, \eps}(T_2^{-A} D^{\frac{5}{8}+\eps})$ for \eqref{error_expression_eq}.
            
            \item[\textit{Sub-case 2c:}] $1 < x < T_2 \log{T_2}$.
            
            Finally, \eqref{L_bound_3} shows that this final sub-case contributes $O(R^\frac{3}{2} T_1^{-2} (D_2/m)^\frac{1}{4})$ to \eqref{error_integral_eq_2}, which translates to $O_\eps(R^\frac{3}{2} T_1^{-2} D^{\frac{5}{8}+\eps})$ in \eqref{error_expression_eq}. This is $O_\eps(D^{\frac{1}{2}} R^{3+\eps})$ by \eqref{reference_eq}, so we have exhausted all possible cases and the proof of \autoref{shifted_conv_lemma} (and therefore also of \autoref{main_thm}) is complete.
            
        \end{itemize}
    \end{itemize}
    
\end{proof}

\section{Limitations and connections to subconvexity}\label{subconvexity_section}

As \autoref{spectral_expansion_lemma} shows and we use in the course of our argument, bounds towards subconvexity have implications to (at least upper bounds for) the variance $\Var(r, R;\Lambda_D)$. We remark that the opposite is also true, in the sense that upper bounds of the correct order of magnitude for the variance imply subconvexity for certain $L$-functions. This clarifies the obstacles for improving \autoref{main_thm}.

For simplicity consider the case of balls, $r = 0$. If one has an upper bound of the (expected) correct order of magnitude for the variance, i.e. $\Var(0, R; \Lambda_D) \ll \sqrt{D} L(1, \chi_D) R^3$, then assuming $R = o(1)$ the argument in \autoref{SH_asymp_lemma} shows that $h_{0, R}(t) \gg R^2$ for $|t|\leq \frac{1}{R}$, so by \autoref{spectral_expansion_lemma} and non-negativity of the terms we get
\begin{equation}\label{moment_implication_eq}
    \sum_{\substack{f \in \mathcal{B}_0(\Gamma) \\ |t_f| \leq \frac{1}{R}}} \frac{L\left(\frac{1}{2}, f\right) L\left(\frac{1}{2}, f \otimes \chi_D\right)}{L\left(1, \sym^2 f\right) |t_f|}+ \frac{1}{2 \pi} \int\limits_{|t| \leq \frac{1}{R}} \left| \frac{\zeta\left(\frac{1}{2} + it\right) L\left( \frac{1}{2} + it, \chi_D \right)}{\zeta(1 + 2it)} \right|^2 \frac{d t}{|t|+1} \ll \frac{L(1, \chi_D)}{R} 
\end{equation}
for squarefree fundamental discriminants $D>0$ (observe that $|t_f| \gg 1$ for $\Gamma$). As an aside, we note that here the significance of the exponent $5/12$ in \autoref{main_thm} becomes clear. This is because the hardest range in \eqref{moment_implication_eq} is $|t_f|, |t| \asymp D^{\frac{1}{12}}$, where the bounds of \autoref{moment_bound_lemma} intersect, and the best one can do is use H{\"o}lder's inequality against the third moment result of \cite{Young} and the large sieve, obtaining
\begin{equation}\label{first_moment_worst_range_eq}
    \sum_{\substack{f \in \mathcal{B}_0(\Gamma) \\ |t_f| \asymp D^\frac{1}{12}}} \frac{L\left(\frac{1}{2}, f\right) L\left(\frac{1}{2}, f \otimes \chi_D\right)}{L\left(1, \sym^2 f\right)}+ \frac{1}{2 \pi} \int\limits_{|t| \asymp D^\frac{1}{12}} \left| \frac{\zeta\left(\frac{1}{2} + it\right) L\left( \frac{1}{2} + it, \chi_D \right)}{\zeta(1 + 2it)} \right|^2 \dd t \ll_\eps D^{\frac{1}{2} + \eps}. 
\end{equation}
An improvement in the first moment bound \eqref{first_moment_worst_range_eq} is essentially equivalent to an extension of the range of $R$ in \autoref{main_thm}.

Going back to our point about subconvexity, dropping all but one term in \eqref{moment_implication_eq} and using the bound $L(1, \sym^2 f) \gg_\eps |t_f|^{-\eps}$ of \cite{HL} we get
\begin{equation*}
    L\left(\frac{1}{2}, f\right) L\left(\frac{1}{2}, f \otimes \chi_D\right) \ll_\eps \frac{D^{\eps} |t_f|^{1+\eps}}{R}
\end{equation*}
for $f \in \mathcal{B}_0(\Gamma)$ with $|t_f| \leq \frac{1}{R}$. The conductor of the product of $L$-functions on the left is $\asymp D^2 |t_f|^4$, so if $f \in \mathcal{B}_0(\Gamma)$ is fixed and $R \gg D^{-\frac{1}{3}+\delta}$ for a given $\delta>0$, we would obtain sub-Weyl subconvexity for $f \otimes \chi_D$ in the twist aspect, which is currently an open problem. 

In conclusion, improving the exponent $5/12$ of \autoref{main_thm} requires a better bound for the first moment \eqref{first_moment_worst_range_eq}, and improving it to anything below $1/3$ seems especially difficult at present, as it implies a challenging case of sub-Weyl subconvexity.


\nocite{*}  
\bibliographystyle{abbrv}
\bibliography{references}

\end{document}